\documentclass[a4paper,11pt]{article}

\usepackage[latin1]{inputenc} %CaractâÃ£res accentuâÃ s - iso 8859-1
\usepackage[T1]{fontenc} %glyphes latines et cesure langues latines
\usepackage{lmodern} %police LatinModern, meilleure qualitâÃ  que 
\usepackage[left=3.2cm,top=3.5cm,right=3.2cm,nohead,foot=1cm]{geometry}
\usepackage[english]{babel} 
\usepackage{amsfonts,latexsym,amssymb} %paquets mathâÃ matiques
\usepackage{mathrsfs,amsmath}
\usepackage{euscript} %ams fonts, euler scripts symbols in math mode
\usepackage[all]{xy}
\usepackage{placeins}%attention compatibilitâÃ  \usepackage[pdftex]{graphicx} avec \DeclareGraphicsExtensions?
\usepackage[pdftex]{graphicx}
\DeclareGraphicsExtensions{.pdf}
\graphicspath{{Figures/}}
\usepackage{xcolor} %couleurs, non prâÃ sent dans MikTex ?
\usepackage[numbers]{natbib} %gestion fine des \cite et \bibitem

\usepackage{verbatim}
\usepackage{tikz}   %TikZ is required for this to work.  Make sure this exists before the next line
\usepackage{3dplot} 
\usepackage{float}

\newlength{\defbaselineskip}
\newcommand{\Proof}{{\noindent {\bf Proof:   }}}
\newcommand{\EndProof}{{\hfill $\Box \qquad $ \endtrivlist}\par \vspace{0.5cm}}%

\renewcommand{\le}{\leqslant}
\renewcommand{\ge}{\geqslant}

\newcommand{\N}{\ensuremath{\mathbb{N}}}
\newcommand{\R}{\ensuremath{\mathbb{R}}}

\newcommand{\Z}{\mathbf{Z}}
\renewcommand{\phi}{\varphi}

\newcommand{\eps}{\varepsilon}

\newcommand{\vecA}{\mathbf{A}}

\newcommand{\vecF}{\mathbf{F}}
\newcommand{\vecG}{\mathbf{G}}
\newcommand{\vecI}{\mathbf{I}}
\newcommand{\vecId}{\mathbf{Id}}
\newcommand{\vecK}{\mathbf{K}}

\newcommand{\vecM}{\mathbf{M}}
\newcommand{\vecN}{\mathbf{N}}
\newcommand{\vecO}{\mathbf{O}}

\newcommand{\vecR}{\mathbf{R}}
\newcommand{\vecS}{\mathbf{S}}
\newcommand{\vecT}{\mathbf{T}}
\newcommand{\vecU}{\mathbf{U}}
\newcommand{\vecV}{\mathbf{V}}

\newcommand{\vecX}{\mathbf{X}}
\newcommand{\vecY}{\mathbf{Y}}
\newcommand{\vecZ}{\mathbf{Z}}

\newcommand{\vece}{\mathbf{e}}
\newcommand{\vecf}{\mathbf{f}}

\newcommand{\vecn}{\mathbf{n}}
\newcommand{\vecp}{\mathbf{p}}
\newcommand{\vecq}{\mathbf{q}}
\newcommand{\vecr}{\mathbf{r}}
\newcommand{\vecs}{\mathbf{s}}
\newcommand{\veck}{\mathbf{k}}
\newcommand{\vecu}{\mathbf{u}}
\newcommand{\vecv}{\mathbf{v}}

\newcommand{\vecx}{\mathbf{x}}
\newcommand{\vecy}{\mathbf{y}}

\newcommand{\vecxi}{\boldsymbol\xi}

\newcommand{\vecsigma}{\boldsymbol\sigma}

\newcommand{\Div}{\mbox{ div }}
\newcommand{\ds}{\displaystyle}
\newcommand{\na}{\nabla}
\newcommand{\Om}{\Omega}

    \newtheorem{theorem}{Theorem}[section]
    \newtheorem{lemma}[theorem]{Lemma}

    \newtheorem{proposition}[theorem]{Proposition}
    \newtheorem{definition}{Definition}[section]

    \newtheorem{remark}{Remark}[section]

\def\div{\hbox{div  }}

\title{Rough wall effect on micro-swimmers }%\thanks{This work has been supported by Direction G\'en\'erale de l'Armement (DGA).}}

\author{David G\'erard-Varet\footnote{Institut de Mathématiques de Jussieu et Université Paris 7, B\^atiment Sophie Germain, 75205 Paris Cedex 13,  France.} , \quad Laetitia Giraldi\footnote{Unité de Mathématiques Pures et Appliquées, UMPA, ENS de Lyon, 46 allée d'Italie, 69364 LYON, France. Supported by Direction G\'en\'erale de l'Armement (DGA).}}
\newcommand{\composante}[1]{\vphantom{A}^{[#1]}}
\newcommand{\x}{\composante{1}}

\renewcommand{\vec}[1]{\overrightarrow{#1}}
\newcommand{\pa}{\partial}

\begin{document}
\maketitle

% ----------------------------------------------------------------------
% fragment pris sur Internet pour calculer l'heure:
% TIME OF DAY - Nelson Beebe - http://www.math.utah.edu/~beebe/
		\newcount\hh
		\newcount\mm
		\mm=\time
		\hh=\time
		\divide\hh by 60
		\divide\mm by 60
		\multiply\mm by 60
		\mm=-\mm
		\advance\mm by \time
		\def\hhmm{\number\hh:\ifnum\mm<10{}0\fi\number\mm}
% fin du fragment pris sur Internet, avec mes remerciements âÃÂ° Nelson Beebe
% ----------------------------------------------------------------------
%brouillon temporaire - version du \today, \hhmm \\

\textbf{Abstract} : We study the effect of a rough wall on the controllability of  micro-swimmers made of several balls linked by thin jacks: the so-called 3-sphere and 4-sphere swimmers.  Our work completes the previous work \cite{AlougesGiraldi12} dedicated to the effect of a flat wall. We show that a controllable swimmer (the 4-sphere swimmer) is not impacted by the roughness. On the contrary, we show that  the roughness changes  the dynamics of the 3-sphere swimmer, so that it can reach any direction almost everywhere.
\\

\noindent{\bf{Keywords}} :
%\vspace{.5mm} \noindent 
Low-Reynolds number swimming; self-propulsion; Three-sphere swimmer; rough wall effect; Lie brackets; control theory; asymptotic expansion.
\tableofcontents

\section{Introduction}
Micro-swimming  is a subject of growing interest,  notably for its biological and medical implications:  one can mention the understanding of reproduction processes,  the description of infection mechanisms, or the conception of micro-propellers for drug delivery in the body. As regards its mathematical modeling and analysis,  the studies  by Taylor \cite{Taylor51}, Lighthill \cite{Lighthill75} and Purcell \cite{Purcell77} have been pioneering contributions to a constantly increasing field:  we refer to the recent work of T. Powers and E. Lauga \cite{PowersLauga09} for an extensive bibliography.
\\

Among the many aspects of micro-swimming, the influence of the environment on swimmers dynamics has been recognized by many  biological studies (see for instance  \cite{BerkeAllison08}, \cite{OrMurray09}, \cite{Rothschild63}, \cite{SmithBlake10}, \cite{SmithGaffney09}, \cite{WinetBernstein84}, \cite{WinetBernstein84a}). One important factor in this dynamics is the presence of confining walls. For example, experiments have shown that some microorganisms, like E. Coli, are attracted to surfaces.    
\\

 The focus of this paper is {\em  the effect  of wall roughness  on micro-swimming}. Such effect has  been already recognized in the context of microfluidics, in connection with  superhydrophobic surfaces (\cite{Bocquet}, \cite{Ledesma-AguilarYeomans13}). Moreover, recent studies have highlighted the role of roughness in the dynamics of passive spherical  particles in a Stokes flow: we refer for instance to the study of S. H. Rad and A. Najafi \cite{RadNajafi10a} or to the one of D. G\'erard-Varet and M. Hillairet \cite{GerardVaretHillairet11}. 
\\

We want here to study the impact of a rough wall on  the displacement of micro-swimmers, at low Reynolds number. Our point of view will be theoretical, more precisely based on control theory. Connection between swimming at low Reynolds number and control theory has been emphasized over the last years  (see \cite{AlougesDeSimone13}, \cite{ChambrionMunnier11a}, \cite{GiraldiMartinon13},  \cite{LoheacMunnier12}, \cite{LoheacScheid11}, \cite{Montgomery02}). We shall ponder here on the recent studies \cite{AlougesDeSimone10} and \cite{AlougesGiraldi12}, dedicated to the controllability analysis of particular Stokesian robots, in the whole space and in the presence of a plane wall respectively. We shall here incorporate roughness at the wall, and focus on two classical models of swimmers:  the 3-sphere swimmer (see \cite{AlougesDeSimone10}, \cite{AlougesDeSimone08}, \cite{AlougesGiraldi12}, \cite{GolestanianAjdari08}) and the 4-sphere swimmer (see \cite{AlougesDeSimone10}, \cite{AlougesGiraldi12}). First, we will show that the controllability of the 4-sphere swimmer  (already true near a flat wall) persists with roughness. Then, we will prove that  the rough wall leads the 3-sphere swimmer to reach any space direction. The underlying mechanism is the symmetry-breaking generated by the roughness.  
\\

The paper is divided into three parts. In Section \ref{Sec:Setting}, we introduce the mathematical model for the fluid-swimmer coupling, and we derive from there an ODE for the  swimmer dynamics. In Section \ref{sec:Neumann-To-Dirichlet_map}, we  show that the force field in this ODE is  analytic with respect to the roughness amplitude and swimmer size and position. Combining this property with the results of \cite{AlougesGiraldi12} yields  the controllability of the 4-sphere swimmer "almost everywhere". Section \ref{sec:DtoN} provides an asymptotic expansion of the Dirichlet-to-Neumann operator,  with respect to the roughness amplitude and swimmer's size. This operator is naturally involved in the expansion of the force fields.  Eventually,  we use this expansion and make it truly explicit in Section \ref{sec:Three_sphere}, in the special case of the 3-sphere swimmer. This allows us to show its controllability.   

\section{Mathematical setting}
\label{Sec:Setting}

In this part, we present our mathematical model for the swimming problem.  

\subsection{Swimmers} 
\label{SubSec:Swimmer}
We carry on the study of specific swimmers that were considered in \cite{AlougesDeSimone10} in $\mathbf{R}^3$ and in \cite{AlougesGiraldi12} in an half plane.  These swimmers consist of $N$ spheres $\displaystyle \cup_{l=1}^{N} B_l$  of radii $a$ connected by $k$  thin jacks which are supposed free of viscous resistance. 
The position of the swimmer is described  by a variable  $\vecp \in \mathbf{R}^3 \times SO(3)$, which gives both the coordinates of one point over the swimmer and the swimmer's orientation.
%To simplify the notation, we assume that the $i$-th arm linked the ball $B_c$ and the ball $B_i$. 
Moreover, the shape variable is denoted by a $k$-tuple $\vecxi$:  its $i$th component $\xi_i$  gives the length of  $i$th arm, that can stretch or elongate through time. 
%A stroke consists in changing the lengths of the jacks in a periodic manner. We also denote by 
%$\overrightarrow{\vecxi_i}$ the {\em unit} vector which defines the direction of the $i$th arm. 
Nevertheless,  the directions of the arms are only modified by global rotation of the swimmer.  Let us stress that all the variables above depend implicitly on time, through the transport and deformation of the swimmer.  
%Finally, we call $\mathcal{S} \subset \mathbf{R}^M$ (for a suitable $M \in \mathbf{N}$) the set of admissible states $(\vecxi,\vecp)$.  Typically, this admissible set will encode geometrical constraints, like a no-contact condition between the swimmer and the confining wall, or the fact that the balls can not touch each other.   $\mathcal{S}$ will be a connected smooth submanifold of $\mathbf{R}^M$. 
\\

Many results of our paper apply to the general class of swimmers just described. Nevertheless,  we will pay a special attention to two examples:
\begin{itemize}
\item {\em The 4-sphere swimmer}.  We consider a regular tetrahedron $(\vecS_1,\vecS_2,\vecS_3,\vecS_4)$ with center $\vecO \in \mathbf{R}^3_+$. The $4$-sphere swimmer consists of four balls linked by four arms of fixed directions $\vec{\vecO\vecS_i}$ which are able to elongate and shrink (in a referential associated to the swimmer).
The four ball cluster is completely described by the list of parameters $(\vecxi,\vecp) = (\xi_1, \dots, \xi_4,\vecx_c,\boldsymbol{\mathcal{R}}) \in  (\sqrt{\frac{3}{2}}a,\infty)^4 \times \mathbf{R}^3 \times SO(3)$. It is known that the 4-sphere swimmer is controllable in $\mathbf{R}^3$ and remains controllable in presence of a plane wall (see  \cite{AlougesDeSimone10}, \cite{AlougesGiraldi12}).
This means that it is able to move to any point and with any orientation under the constraint of being self-propelled, when the surrounding flow is dominated by viscosity (Stokes flow). This swimmer is depicted in Fig. \ref{FourSphereSwimmer}.
\begin{figure}[H]
\begin{center}
\includegraphics[scale=0.6]{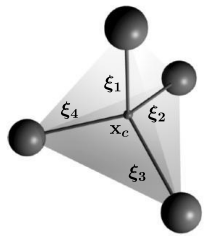}
\end{center}
\caption{The Four-sphere swimmer.\label{FourSphereSwimmer}}
\end{figure}

\item  {\em The 3-sphere swimmer}  (see \cite{AlougesDeSimone10}, \cite{AlougesDeSimone08}, \cite{AlougesGiraldi12} and  \cite{NajafiGolestanian04}). It is composed of three aligned spheres, linked by two arms, see Fig. \ref{Fig:coordinates_intro}. The dynamics of the swimmer is described through the lengths of the two arms $\xi_1, \xi_2$, the  coordinates of the center of the middle ball: $\vecx_c = (x_c,y_c,z_c)$, and some matrix ${\cal R} \in SO(3)$ describing the orientation of the swimmer. Thus, 
$$(\vecxi,\vecp) = (\xi_1,\xi_2, \vecx_c, {\cal R})  \in (2a,\infty)^2 \times \mathbb{R}^3 \times  SO(3).$$  
As regards the position and elongation of the swimmer,  the angle of the rotation  ${\cal  R}$ around the symmetry axis of the 3-sphere is irrelevant. As a matter of fact, we will not show controllability for this angle: our  result, Theorem  \ref{thm:three_sphere}, yields  controllability of the swimmer up to rotation around its axis. Still, the associated angular velocity is not zero, and will appear in the dynamics.   
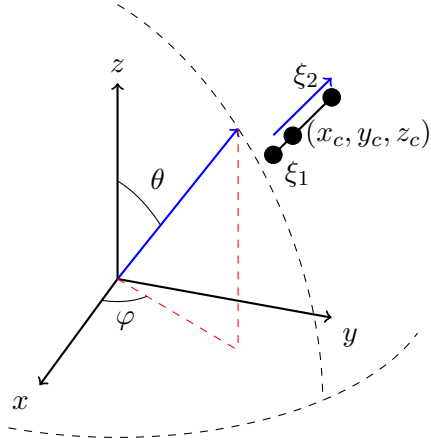
\begin{figure}[H]
\begin{center}
\tdplotsetmaincoords{60}{110}

%define polar coordinates for some vector
%TODO: look into using 3d spherical coordinate system
\pgfmathsetmacro{\rvec}{1.40}
\pgfmathsetmacro{\thetavec}{36}
\pgfmathsetmacro{\thetavecc}{34}
\pgfmathsetmacro{\thetaveccc}{32}
\pgfmathsetmacro{\phivec}{60}
\pgfmathsetmacro{\rvecc}{1.5}
\pgfmathsetmacro{\rveccc}{0.4}
\pgfmathsetmacro{\rvecccc}{0.5}
\pgfmathsetmacro{\rveccccc}{1.1}

%start tikz picture, and use the tdplot_main_coords style to implement the display 
%coordinate transformation provided by 3dplot
\begin{tikzpicture}[scale=3,tdplot_main_coords]

%set up some coordinates 
%-----------------------
\coordinate (O) at (0,0,0);
%\coordinate (V) at (0.9,0.6,0.9);

%determine a coordinate (P) using (r,\theta,\phi) coordinates.  This command
%also determines (Pxy), (Pxz), and (Pyz): the xy-, xz-, and yz-projections
%of the point (P).
%syntax: \tdplotsetcoord{Coordinate name without parentheses}{r}{\theta}{\phi}
\tdplotsetcoord{V}{0.46}{50.33}{34}
\tdplotsetcoord{V}{\rvec}{\thetavec}{\phivec}
\tdplotsetcoord{R}{\rvecc}{\thetavec}{\phivec}
\tdplotsetcoord{S}{\rveccc}{\thetavec}{\phivec}
\tdplotsetcoord{S1}{\rvecccc}{\thetavecc}{\phivec}
\tdplotsetcoord{S2}{\rveccccc}{\thetaveccc}{\phivec}

%draw figure contents
%--------------------

%draw the main coordinate system axes
\draw[thick,->] (0,0,0) -- (1,0,0) node[anchor=north east]{$x$};
\draw[thick,->] (0,0,0) -- (0,1,0) node[anchor=north west]{$y$};
\draw[thick,->] (0,0,0) -- (0,0,1) node[anchor=south]{$z$};

%draw a vector from origin to point (P) 
%\draw[-stealth,color=red] (O) -- (P);
%\draw[dashed,color=red] (O) -- (V);
\draw[thick,->,color=blue] (O) -- (V);

%draw projection on xy plane, and a connecting line
\draw[dashed, color=red] (O) -- (Vxy);
\draw[dashed, color=red] (V) -- (Vxy);

%swimmer
%Monte de 0.1
\draw[fill=black] (0.2,0.8,0.9) circle (0.1em);
\draw[fill=black] (0.5,1,1.2) circle (0.1em);
\draw[fill=black] (1.1,1.4,1.8) circle (0.1em);
\draw[thick] (0.2,0.8,0.9) -- (0.5,1,1.2)--(1.1,1.4,1.8);

\draw  (0.2,0.9,0.9) node [anchor=south west]{$(x_c,y_c,z_c)$};
\draw  (0.3,0.95,0.77) node [anchor=south]{$\xi_1$};
\draw (1.15,1.3,1.8) node [anchor=south]{$\xi_2$};

%Le vecteur
\draw[thick,->,color=blue]  (0.2,0.8,1) -- (1.1,1.4,1.9);

%draw the angle \phi, and label it
%syntax: \tdplotdrawarc[coordinate frame, draw options]{center point}{r}{angle}{label options}{label}
\tdplotdrawarc{(O)}{0.2}{0}{\phivec}{anchor=north}{$\phi$}

%set the rotated coordinate system so the x'-y' plane lies within the
%"theta plane" of the main coordinate system
%syntax: \tdplotsetthetaplanecoords{\phi}
\tdplotsetthetaplanecoords{\phivec}

%draw theta arc and label, using rotated coordinate system
\tdplotdrawarc[tdplot_rotated_coords]{(0,0,0)}{0.5}{0}{\thetavec}{anchor=south west}{$\theta$}

%draw some dashed arcs, demonstrating direct arc drawing
\draw[dashed,tdplot_rotated_coords] (\rvec,0,0) arc (0:90:\rvec);
\draw[dashed] (\rvec,0,0) arc (0:90:\rvec);
\end{tikzpicture}
\end{center}
\caption{Coordinates of the $3$-sphere swimmer}
\label{Fig:coordinates_intro}
\end{figure}

\end{itemize}

%
% and the 4-sphere swimmer (see  \cite{AlougesDesimone10} and \cite{AlougesGiraldi12}).  Here is a description of these swimmers. 
%\\
%
%The Three sphere swimmer is composed of  as shown in Fig. 
%We characterize swimmer's position and orientation by the coordinates $(\theta,\phi,\vecx_c) \in  \mathbf{R}^3\times [0,2\pi]^2$ where $\theta$ is the angle between the swimmer and the $z-$axis, $\phi$ is the angles between the $y$-axis and the projection of the swimmer in a plane $Oxy$ and $\vecx_c$ is the coordinates of the middle spheres. Therefore, in that case, the swimmer is completely described by the vector $\theta,\phi,\vecx_c)\in \mathcal{S}=(2a,\infty)^2 \times \mathbb{R}^3 \times [0,2\pi]^2 \subset (2a,\infty)^2 \times \mathbf{R}^2\times \mathbf{R}/{2\pi\mathbf{Z}}$.
%\\

\subsection{Fluid flow}
We consider a fluid confined by a rough boundary.  This boundary is modelled by a surface with equation $z= \eps h(x,y)$, for some Lipschitz positive function $h$. Here,  $\eps > 0$ denotes the amplitude of the roughness, that is  $\Vert h \Vert_{\infty} = 1$. 
The swimmer evolves in the half-space $\mathcal{O}=\{(x,y,z)\in \mathbf{R}^3\mbox{ s. t. } z>  \eps h(x,y)\}$. The fluid domain is then ${\cal F} \: := \:  \mathcal{O} \setminus \cup_{l=1}^N B_l$, and again it depends implicitly on time. Finally, we assume that the flow is governed there by the Stokes equation.  Thus, the velocity $\vecu^S$ and the pressure $p^S$ of the fluid satisfy:
\begin{equation}
\label{eq:equation_generale}
\begin{array}{ccc}
- \mu \Delta \vecu^S + \nabla p^S = 0\,, & \Div \,\vecu^S = 0 &  \textrm{in ${\cal F}$,}
\end{array}
\end{equation}
where $\mu$ is the viscosity of the fluid. We complement  the Stokes equation \eqref{eq:equation_generale} by standard no-slip boundary conditions, that read: 
\begin{equation}
\label{eq:boundary_condition}
\left\{
\begin{array}{lll}
\vecu^S = \Omega \times (\vecx-\vecx_c) + \vecv + \vecu_d \quad  \textrm{at $\displaystyle\cup_{ l=1}^{N} \partial B_l$,}\\
\vecu^S = 0 \quad \quad  \textrm{at $\partial \mathcal{O}$.}
\end{array}
\right.
\end{equation}
In other words, we impose  the continuity of the velocity both at the fixed wall and at the boundary of the moving swimmer. Note that the velocity field of the swimmer is made of two parts: 
\begin{itemize}
\item one corresponding to an (unknown) rigid movement, with angular velocity $\Omega$ and linear velocity $\vecv$. If $\vecx_c$ is the point over the swimmer  encoded in $\vecp$, the velocity $\vecv$ is its speed. The vector $(\Omega, \vecv)^t$ can be  identified with $\dot{\vecp}$ (everything will be made explicit in due course). 
\item one corresponding to the (known) deformation of the jacks, with associated velocity $\vecu^d$, depending on $\dot{\vecxi}$.   
\end{itemize}

\noindent
Introducing  the Hilbert space
 \begin{equation} \label{hilbert}
 \mathcal{V} =\left\{ \vecu \in \mathcal{D}' (\mathcal{F},\mathbf{R}^3)\, \vert \, \nabla \vecu \in L^2(\mathcal{F}), \, \frac{\vecu(\vecr)}{\sqrt{1+\vert \vecr \vert^2}} \in L^2(\mathcal{F}) \right\},
\end{equation} 
we get (for any configuration of the swimmer $\cup B_l$ and velocities $(\Omega, \vecv, \vecu^d)$) a unique solution $(\vecu^S,p^S)$ of \eqref{eq:equation_generale} -\eqref{eq:boundary_condition} in $\mathcal{V} \times L^2({\cal F})$. See Appendix A for more details.

\subsection{Dynamics}
\label{SubSec:Dynamics}
Of course, the previous relations describe the  equilibrium of the fluid flow at any given instant $t$. To close the model (that is the fluid-swimmer coupling), we still need to specify the dynamics of the swimmer, based on Newton's laws. The description is by now classical  
(see for instance \cite{AlougesDeSimone10},  \cite{LoheacMunnier12}), and can be expressed by  an affine control system without drift. Let us recall the principle of derivation.  Neglecting inertia, Newton's laws become

\begin{equation}
\label{eq:Newton_law}
\left\{
\begin{array}{ll}
\displaystyle \sum_{l=1}^N \int_{\partial B_l} \vecsigma(\vecu^S,p^S)\cdot \vecn \,\mbox{d}s = 0\,, \\
\displaystyle \sum_{l=1}^N \int_{\partial B_l} \vecsigma(\vecu^S,p^S) \cdot \vecn \times (\vecx-\vecx_c) \,\mbox{d}s = 0\,,
\end{array}
\right.
\end{equation}
where $\vecsigma(\vecu,p)= \mu ( \nabla \vecu + \nabla^t \vecu) - p \mathbf{Id}$ is the Cauchy tensor. 
\newline

Moreover, if we introduce an orthonormal basis $(\vece_1,\vece_2,\vece_3)$ and  use linearity, $\vecu^S$ decomposes into
\begin{equation} \label{rigidvf}
\vecu^S = \sum_{i=1}^3\Omega_i \vecu_i + \sum_{i=4}^6 v_{i-3} \vecu_i\, + \vecu^d. 
\end{equation}
Here, the $\vecu_i$'s and $\vecu^d$ are solutions of the Stokes equation, with zero Dirichlet condition at the wall, and inhomogeneous Dirichlet conditions at the ball.  The Dirichlet data is $\vece_i \times (\vecx - \vecx_c)$ for $i=1,2,3$,  $\vece_{i-3}$ for $i=4,5,6$, $\vecu_d$ for $\vecu^d$.
Note also that  the speed $\vecu^d$ can be expressed as a linear combination of $(\dot{\xi}_i)_{i=1}^{k}$: 
\begin{equation}
\label{eq:expression_u_deformation}
\vecu^d = \sum_{i=1}^{k} \vecu^d_i  \dot{\xi}_i.
\end{equation}
%where, for $i \in \left\{1, \cdots,k\right\}$, $\vecu^d_i$ is the Stokes solution associated with the boundary condition, 
%\begin{equation}
%\label{eq:boundary_deformation}
%\left\{
%\begin{array}{lll}
%\vecu^G = 0 \quad  \textrm{in $\displaystyle \cup_{ j=1, j\neq i}^{N} \partial B_j$,}\\
%\vecu^G =  \overrightarrow{\vecxi_i} \,\,\,  \textrm{in $\partial B_i$,}\\
%\vecu^G = 0 \quad  \textrm{in $\partial \mathcal{F}$.}
%\end{array}
%\right.
%\end{equation}
Identifying $(\Omega,\vecv)^t$  with $\dot{\vecp}$ (everything will be made explicit in due course),   the system \eqref{eq:Newton_law} reduces to the following  ODE: 

\begin{equation}
\label{eq:dynamics1}
\vecM(\vecxi,\vecp) \, \dot{\vecp} + \vecN(\vecxi,\vecp) = 0
\end{equation}
where the matrix $\vecM(\vecxi,\vecp)$ is defined by,
$$
\vecM_{i,j}(\vecxi,\vecp):= 
\left\{
\begin{array}{ll}
\displaystyle \sum_{l=1}^N \int_{\partial B_l} \left((\vecx-\vecx_c) \times \vece_i\right) \cdot \vecsigma(\vecu_j,p_j) \,\vecn \,\mbox{d}s \quad \left(1 \leq i \leq 3, 1 \leq j \leq6 \right)\,, \\
\displaystyle \sum_{l=1}^N \int_{\partial B_l} \vece_{i-3} \cdot \vecsigma(\vecu_j,p_j)\, \vecn \, \mbox{d}s \quad \left(4 \leq i \leq 6, 1 \leq j \leq6 \right)\,,
\end{array}
\right.
$$
and $\vecN(\vecxi,\vecp)$ is the vector of $\mathbb{R}^6$ whose entries are,

%$$
%\vecN_i(\vecxi,\vecp):=
%\left\{
%\begin{array}{ll}
%\displaystyle \sum_{l=1, l=c}^{N-1} \sum_{j=1}^{N-1} \int_{\partial B_l} \left(\vecx \times \vece_i\right) \cdot \dot{\vecxi_j} \vecsigma(\vecu^d_j,p^d) \,\vecn_t \,\mbox{d}\vecs \quad \left(1 \leq i \leq 3 \right)\,, \\
%\displaystyle \sum_{l=1, l=c}^{N-1}\sum_{j=1}^{N-1} \int_{\partial B_l} \vece_{i-3} \cdot \dot{\vecxi_j} \vecsigma(\vecu^d_j,p^d) \, \vecn_t \, \mbox{d}\vecs \quad \left(4 \leq i \leq6 \right)\,.
%\end{array}
%\right.
%$$
%
$$
\vecN_i(\vecxi,\vecp):=
\left\{
\begin{array}{ll}
\displaystyle \sum_{l=1}^{N} \int_{\partial B_l} \left((\vecx-\vecx_c) \times \vece_i\right) \cdot\vecsigma(\vecu^d,p^d) \,\vecn \,\mbox{d}s \quad \left(1 \leq i \leq 3 \right)\,, \\
\displaystyle \sum_{l=1}^{N} \int_{\partial B_l} \vece_{i-3} \cdot \vecsigma(\vecu^d,p^d) \, \vecn \, \mbox{d}s \quad \left(4 \leq i \leq6 \right)\,.
\end{array}
\right.
$$

The matrix $\vecM(\vecxi,\vecp)$ is checked to be symmetric and negative definite.  By inverting it in \eqref{eq:dynamics1}, we end up with the following relation for the swimmer's dynamics: 
\begin{equation}
\label{eq:dynamics0}
\dot{\vecp} = - \vecM^{-1}(\vecxi,\vecp) \vecN(\vecxi,\vecp)\,.
\end{equation}
By using \eqref{eq:expression_u_deformation}, we deduce that there are vector fields $\vecF_i$, $i=1..k$, such that the equation \eqref{eq:dynamics0} reads 
\begin{equation}
\label{eq:dynamics}
\dot{\vecp} = \sum_{i=1}^{k} \vecF_i(\vecxi,\vecp) \dot{\xi}_i\,.
\end{equation}

\subsection{Main results}
Before turning to our mathematical analysis, we synthetize here our main results. 
\\

The controllability properties of the swimmers will follow from a careful study of the properties of the $\vecF_i$'s in \eqref{eq:dynamics}. As a first consequence of this study, we will obtain the analyticity of these vector fields with respect to all parameters: the typical height of the roughness $\eps$, the radius of the balls $a$,  the vector of arms lengths $\vecxi$ and  the position of the swimmer $p$.  More precisely, defining 
\begin{multline*}
 {\cal A} \: := \: \{ (\eps, a, \vecxi, \vecp) \in  \R \times  \R_+^* \times  (\R_+^*)^k \times (\R^3 \times SO(3))   \: \mbox{such that} \\ 
 \: B_i \cap  B_j = \emptyset \:  \forall i\neq j, \: \mbox{and} \: B_i \cap \partial \mathcal{O} = \emptyset \: \forall i \},   $$
\end{multline*}
we have the following
\begin{theorem}
\label{thm:analyticity}
For all $i=1\dots k$, the field $\vecF^i(\vecxi,p)$ (which depends also implicitly on $\eps$ and $a$) is an analytic function of $(\eps, a,\vecxi,\vecp)$ over ${\cal A}$.  
\end{theorem}
Then, as a consequence of Theorem \ref{thm:analyticity}, we will prove that the roughness does not change the controllability of the $4$-sphere swimmer. We restrict here to local controllability "almost everywhere":  this terminology refers to the following 
\begin{definition}  {\bf ("almost everywhere")}
We say that a property holds for almost every $(\eps,a,\vecxi,\vecp)$ in ${\cal A}$ if it holds for all  $(\eps,a,\vecxi,\vecp)$  outside the zero set of a (non-trivial) analytic function over ${\cal A}$.
\end{definition}
%{\red Note that for any given values of $\eps$ and $a$,   the admissible set ${\cal S} = {\cal S}_{\eps,a}$ of the swimmer is simply %the set of couples $(\vecxi,\vecp)$ such that $(\eps,a,\vecxi,\vecp)$ belongs to ${\cal A}$. } 
We have
\begin{theorem}
\label{thm:four_sphere}
The $4$-sphere swimmer is controllable almost everywhere, in the following sense:  for almost every  $(\eps,a,\vecxi^i,\vecp^i)$, one has local controllability from the initial configuration $(\vecxi^i,\vecp^i)$. This means that  for any final configuration $(\vecxi^f,\vecp^f)$ in a small enough neighborhood of $(\vecxi^i,\vecp^i)$ and any final time $T>0$, there exists a stroke $\vecxi\in \mathcal{W}^{1,\infty}([0,T])$, satisfying
$\vecxi(0)=\vecxi^i$ and $\vecxi(T)=\vecxi^f$ and such that if the self-propelled swimmer starts in position $\vecp^i$ with the shape $\vecxi^i$ at time $t=0$, it ends at position $\vecp^f$ and shape $\vecxi^f$ at time $t=T$ by changing its shape along $\vecxi(t)$.
\end{theorem}

In the last Section \ref{sec:Three_sphere}, we shall address the controllability of the 3-sphere swimmer. In the case of a flat boundary, as shown in  \cite{AlougesGiraldi12}, symmetries constrain the swimmer to move in a plane. Also, it does not rotate around its own axis. As we will see, the roughness at the wall breaks (in general) such symmetries, allowing for local controllability almost everywhere. Let us point here a subtlety regarding our controllability result. To express the dynamics of the swimmer through the equation \eqref{eq:dynamics}, we have included in variable $\vecp$ (more precisely in its $SO(3)$ component)  an angle describing rotation of the $3$-sphere around its own axis. We are not able to show controllability for this angle: we only show controllability for the other components of $\vecp$. Of course, this is not a problem with regards to the effective movement of the swimmer:  this angle is indeed irrelevant with regards to the swimmer's orientation and position.  The analysis of Section \ref{sec:Three_sphere} leads to the  
 \begin{theorem}
 \label{thm:three_sphere}
There exists a surface $h \in C^\infty_c(\R^2)$ such that the 3-sphere swimmer is locally controllable almost everywhere (up to rotation around its axis). 
% Assume (for instance) that the rough surface is described by   $h(x,y) = $ \textbf{mettre ici la fonction utilisee dans MAPLE}. Then,   
 \end{theorem}
Refined statements will be provided in Section \ref{sec:Three_sphere}. This controllability result requires a careful asymptotic asymptotic expansion of the force fields $\vecF^i$. This expansion is related to an expansion of a Dirichlet-to-Neumann map, performed in section \ref{sec:DtoN}. Eventually, the dimension of the Lie algebra generated by the force fields is computed numerically, and the controllability result follows from application of Chow's theorem.

\section{Analyticity of the dynamics}
\label{sec:Neumann-To-Dirichlet_map}

\subsection{Regularity}
This paragraph is devoted to the proof of Theorem \ref{thm:analyticity}. Let $\overline{\vecY} = (\overline{\eps}, \overline{a}, \overline{\vecxi}, \overline{\vecp}) \in {\cal A}$.  We must prove analyticity of the $\vecF^i$'s with respect to ${\vecY} = (\eps, a, \vecxi, \vecp)$, in a neighborhood of $\overline{\vecY}$. It will follow from the analyticity of $\vecM$ and $\vecN$ defined after \eqref{eq:dynamics1}. Their definitions involve  functionals of the type 
$$ \vecI \: := \:  \sum_{l=1}^N \int_{\pa B_l} \left( \begin{smallmatrix} 1\\ \vecx \end{smallmatrix}\right)   \otimes \sigma(\vecu,p) n \, ds $$
 where $(\vecu,p)$ satisfies the Stokes equation in ${\cal F}$, with Dirichlet conditions of the type: 
 $$ \vecu = 0 \quad  \mbox{at} \: \pa {\cal O}, \quad  \vecu = \vecu_l \quad \mbox{at} \: \pa B_l, \quad l=1,\dots,N $$
 for some family of rigid fields  $\vecu_l$'s  taken in the "elementary set" 
 $\{\vece_i \times \vecx, \vece_i, \: i=1...3 \}$.      
 \\
 
We denote by $\vecx_l$, resp. $\overline{\vecx}_l$  the center of the ball $B_l$, resp. the center of the ball $\overline{B}_l$ corresponding to $\overline{\vecY}$.   We  introduce the diffeomorphisms 
$$\varphi_l(\vecx) \: := \:  \frac{a}{\overline{a}} (\vecx - \overline{\vecx}_l) \: + \: \vecx_l. $$
Then, we have 
$$  \int_{\pa B_l} \left( \begin{smallmatrix} 1\\ \vecx \end{smallmatrix}\right)   \otimes \sigma(\vecu,p) n \, ds \: = \: \left( \frac{a}{\overline{a}} \right)^2 \int_{\pa \overline{B}_l}   \left( \begin{smallmatrix} 1\\ \varphi_l(\vecx) \end{smallmatrix}\right)   \otimes \sigma(\vecu \circ \varphi_l,p \circ \varphi_l) n \, ds. $$
Hence,  in order to prove Theorem \ref{thm:analyticity}, it is enough to show that for all $l=1...N$, for  $\delta, \eta > 0$ small enough:
$$ B(\overline{\vecY}, \delta) \:  \mapsto  \: H^1\left({\cal F} \cap B(\overline{\vecx}_l, a+\eta)\right) \times L^2\left({\cal F} \cap B(\overline{\vecx}_l, a+\eta)\right), \quad \vecY \mapsto (\vecu \circ \varphi_l, p \circ \varphi_l) $$
is analytic. Indeed, $\vecY \mapsto \sigma(\vecu \circ \varphi_l,p \circ \varphi_l)$ will  be analytic with values in $H^{-1/2}(\pa \overline{B}_l)$, and the surface integral will be analytic as well. 
\\

Therefore, we define the change of variable 
\begin{equation*}
\varphi(\vecx) = \vecx +  \sum_l \chi(\vecx - \overline{\vecx}_l) \left( \varphi_l(\vecx) - \vecx \right)    + (\eps - \overline{\eps})  \chi_h(\vecx)  (0,0,h(x_1,x_2)) 
\end{equation*}
with  $\chi, \chi_h  \in C^\infty_c(\R^3)$, $\chi = 1$ near $B(0,\overline{a})$, $\chi_h = 1$ near $x_3 = \overline{\eps} \, h(x_1,x_2)$. 
For $\chi$ and $\chi_h$ with small enough supports, and for $\displaystyle \vecY \in B(\overline{\vecY}, \delta) $, $\delta > 0$ small enough, it is easily seen that $\varphi$ is a smooth diffeomorphism, which depends analytically on $\vecY$, and such that 
$\displaystyle  \varphi(\overline{{\cal F}}) \: = \: {\cal F}.$ Moreover, one has $\varphi = \varphi_l$  in a small enough $\delta'$-neighborhood of $\overline{B}_l$.  Introducing $\vecU := \vecu \circ \varphi$ and $P := p \circ \varphi$, it remains to prove the following 
\\

{\em Claim:} $\vecY \mapsto \vecU$ is analytic from $B(\bar{\vecY}, \delta)$ to $\overline{{\cal V}}_0$, where 
$$\overline{\mathcal{V}}_0 \: := \: \left\{ \vecU \in \mathcal{D}' (\overline{\mathcal{F}},\mathbf{R}^3)\, \vert \, \nabla \vecU \in L^2(\overline{\mathcal{F}}), \, \frac{\vecU(\vecr)}{\sqrt{1+\vert \vecr \vert^2}} \in L^2(\overline{\mathcal{F}}), \quad \vecU\vert_{\pa \bar{{\cal O}}} = 0\right\}. $$
To prove this claim, one first needs to write down the system satisfied by $\vecU, P$. A simple computation yields
\begin{equation} \label{stokesAB}
\left\{
\begin{aligned}
-\Div(A \na \vecU) + B \na P  & = 0 \quad \mbox{in} \: \overline{{\cal F}}, \\
\Div (B^t \vecU) & = 0 \quad \mbox{in} \: \overline{{\cal F}},  \\
\vecU = 0   \quad \mbox{at} \: \pa \overline{{\cal O}}, \quad \vecU & = \vecU_l \quad \mbox{at} \: \pa \overline{{\cal B}}_l, 
\end{aligned}
\right.
\end{equation}
where 
\begin{align*}
& A = A(\vecx) := |\det \na \phi(\vecx) | (\na \phi^{-1})^t (\na \phi^{-1})(\phi(\vecx)), 
\\ 
& B = B(\vecx) :=   |\det \na \phi(\vecx) | (\na \phi^{-1})(\phi(\vecx)), \quad \vecU_l := \vecu_l \circ \varphi_l.  
\end{align*}
Note that $A,B, \vecU_l$ depend analytically on the parameter $\vecY$. We now introduce 
$$ \overline{\mathcal{V}}' \: := \: \textrm{the dual space of} \quad \left\{ \vecU  \in \overline{\mathcal{V}}_0, \quad \vecU\vert_{\pa B_l} = 0, \quad  l=1 \dots N  \right\}, $$ 
and consider the mapping 
\begin{align*}
&{\cal L} : \:  \: B(\overline{\vecY}, \delta) \times \overline{{\cal V}}_0 \times L^2(\overline{\cal F}) \: \mapsto \:  \overline{{\cal V}}' \times L^2(\overline{\cal F}) \times \prod_{l} H^{1/2}(\pa \overline{B}_l), \\
& \left(\vecY, \vecV, Q\right) \mapsto \bigl( -\Div (A \na \vecV) + B \na Q, \Div(B^t \vecV), \bigl( \vecV\vert_{\pa \overline{\cal B}_l} - \vecU_l \bigr)_{l=1}^N\bigr). 
\end{align*} 
${\cal L}$ is clearly well-defined, and it is analytic in $(\vecY, \vecV,Q)$: we refer to \cite{Whittlesey} for the definition of analytic functions over Banach spaces. Moreover,  $\vecU = \vecU_{\vecY}$ and $P = P_\vecY$ satisfy 
$$ {\cal L}(\vecY, \vecU, P) = 0 $$
By the analytic version of  the implicit function theorem, see again \cite{Whittlesey},  $\vecU$ and $P$ will be analytic  in $\vecY$ near $\overline{\vecY}$ if 
$$ \frac{\pa {\cal L}}{\pa (\vecV,Q)}\vert_{(\overline{\vecY}, \vecU, P)} \: \mbox{is an isomorphism from} 
\: \overline{{\cal V}}_0 \times L^2(\overline{\cal F}) \: \mbox{to} \:  \overline{{\cal V}}' \times L^2(\overline{\cal F}) \times \prod_{l} H^{1/2}(\pa \overline{B}_l).$$
In other words, analyticity of $\vecU$ and $P$ follows from the existence and uniqueness in 
$\overline{{\cal V}}_0 \times L^2(\overline{\cal F})$ of a solution $(\vecV,Q)$  for the Stokes system 
\begin{align*}
-\Delta \vecV  + \na Q & = \vecF \quad  \: \mbox{in} \: \overline{{\cal F}}, \\
 \Div \vecV & =  G \quad  \: \mbox{in} \: \overline{{\cal F}}, \\ 
 \vecV & = \vecV_l   \quad \mbox{at} \: \pa \overline{B}_l, \quad l=1...N 
\end{align*}
where $\vecF \in  \overline{{\cal V}}'$, $G \in L^2(\overline{\cal F})$ and $\vecV_l \in H^{1/2}(\pa \overline{B}_l)$ are prescribed data. Note that the space $\overline{{\cal V}}_0$  encodes the additional boundary condition: $\vecV = 0$ at $\pa \overline{{\cal O}}$.    

\medskip
The  well-posedness of the previous system is established in the appendix. This ends the proof of Theorem \ref{thm:analyticity}.  

\subsection{Application to the 4-sphere swimmer}

%\david{A RETOUCHER UN PEU}

%\begin{theorem}
%\label{thm:four_sphere}
%Consider the Four-sphere swimmer depicted in Fig \ref{FourSphereSwimmer} and described in Section \ref{SubSec:Swimmer}, and assume it is self-propelled in a three dimensional viscous flow modeled by Stokes equations in the half space $\mathcal{F}$. Then for almost any initial configuration $(\vecxi^i,\vecp^i)\in \mathcal{S}$, any final configuration $(\vecxi^f,\vecp^f)$ in a suitable neighborhood of $(\vecxi^i,\vecp^i)$ and any final time $T>0$, there exists a stroke $\vecxi\in \mathcal{W}^{1,\infty}([0,T])$, satisfying
%$\vecxi(0)=\vecxi^i$ and $\vecxi(T)=\vecxi^f$ and such that if the self-propelled swimmer starts in position $\vecp^i$ with the shape $\vecxi^i$ at time $t=0$, it ends at position $\vecp^f$ and shape $\vecxi^f$ at time $t=T$ by changing its shape along $\vecxi(t)$.
%\end{theorem}
From the analyticity shown above and the results of  \cite{AlougesGiraldi12}, we can deduce Theorem \ref{thm:four_sphere}. First, by \eqref{eq:dynamics}, we can write the swimmer's dynamics as 
$$ \dot{\left( \begin{smallmatrix} \vecxi \\ \vecp \end{smallmatrix} \right)} = \sum_{i=1}^4 \vecG_i\left(\left( \begin{smallmatrix} \vecxi \\ \vecp \end{smallmatrix} \right)\right) u_i  $$ 
where $(u_i := \dot{\xi}_i)_{i=1}^4$ is the family of controls, and $\vecG_i := \left( \begin{smallmatrix} \vece_i \\ \vecF_i \end{smallmatrix} \right)$ ($(\vece_1,...,\vece_4)$ is the canonical basis of $\R^4$).   By the analyticity of the $\vecG_i$'s and Chow's theorem, it is then enough to prove that for some $(\eps, a, \vecxi,\vecp) \in {\cal A}$,  
$$ \dim Lie_{(\vecxi,\vecp)}(\vecG_1,...,\vecG_{10}) = 10.  $$
We write
$$ \pa^\alpha \vecG_i(\vecxi,\vecp) = \pa^\alpha \vecG_i^0(\vecxi,\vecp) \: + \: O(\eps), \quad \forall \alpha \in \N^7,
$$ 
where the $\vecG_i^0$'s are force fields corresponding to the flat case $h=0$. In particular, for $\eps$ small enough
$$ \dim Lie_{(\vecxi,\vecp)}(\vecG_1,...,\vecG_{10}) \: \ge \:  \dim Lie_{(\vecxi,\vecp)}(\vecG^0_1,...,\vecG^0_{10}).$$
But from \cite{AlougesGiraldi12} we know that for almost every $(a,\vecxi,\vecp)$ 
$$  \dim Lie_{(\vecxi,\vecp)}(\vecG^0_1,...,\vecG^0_{10}) = 10. $$ 
This concludes the proof. 

\section{Asymptotic expansion of the Dirichlet-to-Neumann} \label{sec:DtoN}
We now turn to  the controllability properties of the 3-sphere swimmer. As before, the key point is to determine the dimension of the Lie algebra generated by the force fields $\vecF^i$. Therefore, we need to derive an asymptotic expansion of the $\vecF^i$'s, in $a$ and $\eps$. \\

{\em A preliminary step is to derive an asymptotic expansion of the so-called Dirichlet-to-Neumann map of the Stokes operator.}      Indeed,  the force fields $\vecF^i$ involve this map:  that is, the definition of the coefficients $\vecM_{ij}$ and $\vecN_i$ involves 
$$ DN: \: \prod_{l=1}^N  H^{1/2}(\pa B_l) \mapsto \prod_{l=1}^N  H^{-1/2}(\pa B_l), \quad (\vecu_l) \mapsto \left(\vecf_l \: := \: \sigma(\vecu,p)n\vert_{\pa B_l} \right),  $$
  where $(\vecu,p)$ is the solution of the Stokes equation 
  $$ -\Delta \vecu + \na p = 0, \quad \Div \vecu = 0  \quad \textrm{in} \: {\cal F}, \quad \vecu\vert_{\pa {\cal O}} = 0, \quad \vecu\vert_{\pa B_l} = \vecu_l. $$
More precisely, it involves $DN$ {\em in restriction to $N$-uplets of  rigid  vector fields over $B_l$}, $l=1...N$. We denote by $R$ the (finite-dimensional) space of such $N$-uplets.  
\\

Even restricted to $R$, this operator is not very explicit: to derive directly an expansion in terms of the parameters of the swimmer and wall is not easy. Hence, we follow the same path as in \cite{AlougesDeSimone10,AlougesGiraldi12}: we write that for all  $(\vecu_l)_{l=1}^N \in R$, 
$$    DN\left((\vecu_l)\right)  \: = \: T^{-1}\left( (\vecu_l)\right)$$ 
where 
$$T : \prod_{l=1}^N  H^{-1/2}(\pa B_l) \mapsto  \prod_{l=1}^N  H^{1/2}(\pa B_l), \quad (\vecf_l) \mapsto (\vecu_l \: := \: \vecu\vert_{\pa B_l})  $$
and $\vecu$ is the solution of the following Stokes system in ${\cal O}$:
$$ -\Delta \vecu + \na p \: = \:  \sum_{l=1}^N 1_{\pa B_l}  \vecf_l, \quad \Div \vecu = 0  \quad \mbox{in} \quad {\cal O}, \quad \vecu\vert_{\pa {\cal O}} = 0.$$ 
Equivalently, this last system can be written: 
$$ -\Delta \vecu + \na p = 0, \quad \Div \vecu = 0 \quad \mbox{in} \quad   {\cal O}\setminus \cup_l \pa B_l,  \quad [\vecu]\vert_{\pa B_l} = 0, \quad [\sigma(\vecu,p)n]\vert_{\pa B_l} = \vecf_l, $$   
where $[ \, ]\vert_{\pa B_l}$ denotes the jump across $\pa B_l$. 
Let us remind that ${\cal O} = \{ z > \eps h(x,y) \}$ is the domain without the balls. In particular, the operator $T$ (associated to a transmission condition) is not the Neumann-to-Dirichlet operator. The latter one would correspond to the Stokes problem 
$$ -\Delta \vecu + \na p = 0, \quad \Div \vecu = 0 \quad \mbox{in} \quad   {\cal O}\setminus \cup_l B_l,  \quad \sigma(\vecu,p)n\vert_{\pa B_l} = \vecf_l, $$ 
associated to a Neumann type condition.  However, {\em in restriction to the space $R$, the operators $DN$ and $T^{-1}$ coincide}, due to the fact that a rigid vector field is a solution of the Stokes equation, with zero pressure and zero stress tensor. 
\\

The advantage of  $T$ over the Neumann-to-Dirichlet operator is its more explicit representation. Indeed, one has  for all $i = 1...N$ 
$$  T(\vecf)_i(\vecx) \: = \: \sum_{l=1}^n \int_{\pa B_l} \vecK^\eps(\vecx,\vecy) \vecf_l(\vecy) d\vecy, \quad \vecx \in \pa B_i, $$  
where the kernel $\vecK^\eps$ is simply the Green function associated to the Stokes equation in ${\cal O}$: in other words, 
$(\vecK^{\epsilon},\vecq^{\epsilon})$ is the solution of the problem: 
\begin{equation}
\label{eq:KEquationDirac}
\left\{
\begin{aligned}
- \mu \Delta_{\vecx} \vecK^{\epsilon}(\vecx,\vecx_0) + \nabla_{\vecx} \vecq^{\epsilon}(\vecx)  & =   \delta_{\vecx_0}(\vecx) \,  \vecI, \quad  \vecx  \: \textrm{in} \:  \mathcal{O}, \\ 
{\rm div}_{\vecx} \,\vecK^{\eps}(\vecx,\vecx_0) & = 0, \quad \vecx \:  \textrm{in} \:  \mathcal{O},\\
\vecK^{\epsilon}(\vecx,\vecx_0) &  = 0  \quad \vecx \: \textrm{on} \:  \partial \mathcal{O},
\end{aligned}
\right.
\end{equation} 
where $\vecI$ stands for the identity matrix. This will make easier the derivation of an asymptotic expansion, through an expansion of $T$. Still, there is one little technical difficulty: the domain of definition and range of $T$, that are $\prod_l H^{\pm 1/2}(\pa B_l)$ depend on the parameter $a$ (and also on $(\vecp,\vecxi))$.  Let us denote $B : = B(0,1)$ the unit ball, and $H_N^{\pm1/2} \: := \: \left( H^{\pm 1/2}(\pa B) \right)^N$. We introduce 
$$ \phi : \prod_{l=1}^N H^{1/2}(\pa B_l) \rightarrow  H^{1/2}_N, \quad  \vecu = (\vecu_l) \mapsto \vecU = (\vecU_l : \vecr \mapsto  \vecu_l(\vecx_l + a \vecr)), $$
 as well as the adjoint map 
 $$ \phi^* : H^{-1/2}_N \rightarrow  \prod_{l=1}^N H^{-1/2}(\pa B_l), \quad   \vecF = (\vecF_l) \mapsto \vecf = (\vecf_l), $$
defined through the duality relation: $<\phi^*(\vecF), \vecu>  \: = \:   <\vecF, \phi(\vecu)>$.
Finally, we set  ${\cal T} \: := \:  \phi \circ T \circ \phi^* :  H^{-1/2}_N \mapsto H^{1/2}_N$. We shall use ${\cal T}$ rather than $T$ to compute the expansion of the force field in section \ref{Subsec:Self-propulsion}. Note that ${\cal T}$ depends implicitly on $\eps$, $a$ and  on $(\vecp,\vecxi)$. {\em In what follows, we will always consider configurations in which the swimmer stays away from the rough wall:  
\begin{equation} \label{stayaway}
 dist(B_l,\pa {\cal O})  \ge \delta  > 0, \quad \forall l=1...N,
 \end{equation}
for some given $\delta$.}

\subsection{Expansion for small  $\eps$}
Under the constraint \eqref{stayaway}, we prove
\begin{proposition}  \label{prop:Green_function}
$$ {\cal T} \: :=  \:  {\cal T}^0 \: + \: \eps {\cal T}^1  \: + \: O(\eps^2) \quad  \mbox{in} \quad {\cal L}(H^{-1/2}_N, H^{1/2}_N)  $$ 
where ${\cal T}^0$ and ${\cal T}^1$ are defined in \eqref{Tau0} and \eqref{Tau1}-\eqref{K1} respectively. 
\end{proposition}
{\em Proof.} For $\vecf = (\vecf_l) \in H^{-1/2}_N$, we can write 
$$ {\cal T}(\vecf)_i(r)  \: = \: \sum_j  \int_{\pa B} \vecK^\eps(x_i + a \vecr, x_j + a \vecs) \vecf_j(\vecs) d\vecs $$
(with a classical and slightly abusive notation: the integral should be understood as a duality bracket). Thus, the whole point is to expand the kernel $\vecK^\eps$ defined in \eqref{eq:KEquationDirac}. Of course, the first term should be $\vecK^0$, that is the Green function in the flat case. This Green function can be computed in terms of the Stokeslet by the method of images (see \cite{Blake71}): one has 
\begin{equation}
\label{DefK}
\vecK^0(\vecr,\vecr_0)=\vecG(\vecr-\vecr_0)+\vecK_1(\vecr,\vecr_0)+\vecK_2(\vecr,\vecr_0)+\vecK_3(\vecr,\vecr_0)\,,
\end{equation}
the four functions $\vecG$, $\vecK_1$, $\vecK_2$ and $\vecK_3$ being respectively
the Stokeslet
\begin{equation}
\label{stokeslet}
\vecG(\vecr) = \frac{1}{8\pi \mu} \left( \frac{\vecId }{| \vecr  |} + \frac{\vecr \otimes \vecr}{\vert \vecr \vert^3} \right)
\end{equation}
and the three ``images''
\begin{equation}
\label{DefK1}
\vecK_1(\vecr,\vecr_0) = -\ds \frac{1}{8 \pi \mu} \left( \frac{\vecId}{\vert \vecr' \vert} + \frac{\vecr'\otimes \vecr'}{\vert \vecr' \vert^3} \right)\,,
\end{equation}
\begin{equation}
\label{DefK2}
K_{2,ij}(\vecr,\vecr_0) = \frac{1}{4\pi \mu}  z_0^2 \left(1 - 2\delta_{j3}\right) \left( \frac{\delta_{ij}}{\vert \vecr' \vert^3} - \frac{3 r'_i r'_j}{\vert \vecr' \vert^5} \right)\,,
\end{equation}
\begin{equation}
\label{DefK3}
K_{3,ij}(\vecr,\vecr_0) = -\frac{1}{4\pi \mu}  z_0 \left(1 - 2 \delta_{j3} \right) \left( \frac{r'_3}{\vert \vecr' \vert^3}\delta_{ij} - \frac{r'_j}{\vert \vecr' \vert^3}\delta_{i3} + \frac{r'_i}{\vert \vecr' \vert^3} \delta_{j3} - \frac{3 r'_i r'_j r'_3}{\vert \vecr' \vert^5} \right) \,.
\end{equation}
Here $\vecr_0=(x_0,y_0,z_0)$ and $\vecr' = \vecr - \tilde{\vecr}_0$, where $\tilde{\vecr}_0=(x_0,y_0,-z_0)$ stands for the ``image'' of $\vecr_0$, that is to say, the point symmetric to $\vecr_0$ with respect to the flat wall. \\

We now consider $\vecu^\eps(\vecx, \vecx_0) = \vecK^\eps(\vecx,\vecx_0) - \vecK^0(\vecx,\vecx_0)$, for $\vecx_0 \in \cup_l B_l$. As a function of $\vecx$, it satisfies the Stokes equation in ${\cal O}$:
$$ -\Delta \vecu^\eps(\cdot, \vecx_0) + \na \vecp(\cdot, \vecx_0) = 0, \quad \Div \vecu^\eps(\cdot, \vecx_0) = 0 \quad  \mbox{in} \: {\cal O} $$
with Dirichlet condition 
$$ \vecu^\eps(\cdot,\vecx_0) = - \vecK^0(\cdot, \vecx_0),  \quad  \mbox{at} \quad  \pa {\cal O}.  
$$
We can then expand the boundary data: for $\vecx =  (x,y, \eps h(x,y)) \in {\cal O}$
$$ -\vecK^0(\vecx,\vecx_0) \: = \: - \sum_{k=1}^n \eps^k \frac{h(x,y)^k}{k!} \pa_z^k \vecK^0(x,y,0,\vecx_0) \: + \: O(\eps^{n+1}). $$   
More precisely, under the constraint \eqref{stayaway}, one has 
$$ \| -\vecK^0(\cdot, x_0)  + \sum_{k=1}^n \eps^k \left( \vecx \mapsto \frac{h(x,y)^k}{k!} \pa_z^k \vecK^0(x,y,0,\vecx_0) \right) \|_{H^s(\pa {\cal O})} \: \le \: C_{\delta, s} \,  \eps^{n+1}, \quad \forall \, s. $$
We deduce  from this inequality that 
\begin{equation}  \label{estimetracen}
\| \na \bigl(  \vecu^\eps(\cdot, \vecx_0) - \sum_{k=1}^n  \eps^k \vecu^k(\cdot,\vecx_0) \bigr) \|_{L^2({\cal O})} \: \le \:   C \eps^{n+1}  
\end{equation}
where $\vecu^k$ is the solution of 
\begin{multline*}
 -\Delta \vecu^k(\cdot,\vecx_0) + \na \vecp(\cdot,\vecx_0) = 0, \quad \Div \vecu^k(\cdot,\vecx_0) = 0 \quad  \mbox{in} \quad  {\cal O}, \\
  \vecu^k(\vecx,\vecx_0) =  -  \frac{h(x,y)^k}{k!} \pa_z^k \vecK^0(x,y,0,\vecx_0), \quad \vecx \in \pa {\cal O}.   
  \end{multline*}
The existence of the $\vecu^k$'s and the estimate \eqref{estimetracen}  are obtained by classical arguments (see the appendix for the more difficult case of a rough half-space minus the balls). In particular, we have 
 \begin{equation}  \label{estimetrace1}
\| \na \bigl(  \vecu^\eps(\cdot, \vecx_0) - \eps \vecu^1(\cdot,\vecx_0) \bigr) \|_{L^2({\cal O})} \: \le \:   C \eps^{2}.
\end{equation}
The last step consists in replacing $\vecu^1$ by the solution $\vecK^1$ of 
\begin{multline*}
 -\Delta \vecK^1(\cdot,\vecx_0) + \na \vecp(\cdot,\vecx_0) = 0, \quad \Div \vecK^1(\cdot,\vecx_0) = 0, \quad z > 0, \\
  \vecK^1(x,y,0,\vecx_0) =  -  h(x,y) \pa_z \vecK^0(x,y,0,\vecx_0), \quad (x,y) \in \mathbb{R}^2 , 
  \end{multline*}
that is replacing the rough half-space by the flat half-space. We claim that 
$$ ||\na(\vecu^1(\cdot, \vecx_0) - \vecK^1(\cdot, \vecx_0))\|_{L^2({\cal O} \cap \{ z > 0 \})} \: = \: O(\eps^2).  $$
With no loss of generality, we can assume that $h > 0$ (meaning that the flat wall is below the rough wall). Otherwise, we can make an intermediate comparison with the solution $\tilde\vecK^1$ of the same Stokes problem  in $\{ z > -\eps (\sup |h| + 1) \}$. Now, an easy but important remark is that   
$$ \| \vecK^1(\cdot, \vecx_0)  \|_{H^s(\{0 < z < Z \})} \: \le \: C_{s,Z},  \quad \forall s \in \N, \: \forall Z >0.  $$
 Hence, 
 $$ \vecK^1(\vecx, \vecx_0) = -  h(x,y) \pa_z \vecK^0(x,y,0,\vecx_0) \: + \: O(\eps)  \quad \mbox{in} \quad  H^s(\pa {\cal O}). $$
 By a simple estimate on $\vecu^1 - \vecK^1$, we deduce the claim. 
 \\
 
 Back to the definition of $\vecu^\eps$, we obtain  thanks to standard elliptic regularity in variable $\vecx$: for all $\alpha \in \N^3$,  
 $$ | \pa^\alpha_{\vecx} \left(\vecK^\eps(\vecx,\vecx_0) -  \vecK^0(\vecx,\vecx_0) - \eps \vecK^1(\vecx,\vecx_0) \right)  |  = O(\eps^2),$$ 
  uniformly  in $\vecx, \vecx_0 \in \cup_l B_l$. The same reasoning as above can then be applied to the fields $\vecu^\eps_\beta = \pa^\beta_{\vecx_0} (\vecK^\eps - \vecK^0)$, for all $\beta \in \N^3$.  Hence,  
 $$ |  \pa^\alpha_{\vecx}   \pa^\beta_{\vecx_0} \left( \vecK^\eps(\vecx,\vecx_0) -  \vecK^0(\vecx,\vecx_0) - \eps \vecK^1(\vecx,\vecx_0) \right)|  = O(\eps^2),$$ 
  uniformly  in $\vecx, \vecx_0 \in \cup_l B_l$. The theorem follows straightforwardly, considering 
  \begin{equation} \label{Tau0}
   {\cal T}^0(\vecf)_i(r)  \: := \: \sum_j  \int_{\pa B} \vecK^0(x_i + a \vecr, x_j + a \vecs) \vecf_j(\vecs) d\vecs 
   \end{equation}
   and 
     \begin{equation} \label{Tau1}
   {\cal T}^1(\vecf)_i(r)  \: := \: \sum_j  \int_{\pa B} \vecK^1(x_i + a \vecr, x_j + a \vecs) \vecf_j(\vecs) d\vecs. 
   \end{equation}
Expressing  $\vecK^1(\vecx,\vecx_0)$ with a Poisson kernel yields
\begin{equation} \label{K1}
 \vecK^1(\vecx, \vecx_0) \: := \:  - \int_{\partial \R_+^3}\, h(\vecs) \, \frac{\partial}{\partial z}\left(\vecs \mapsto \vecK^0(\vecs,\vecx) \right) \, \frac{\partial}{\partial z}\left(\vecs \mapsto \vecK^0(\vecs,\vecx_0) \right) \mbox{d}\vecs\,.
\end{equation}
where for simplicity we write $h(\vecs)$ instead of $h(s_1,s_2)$,  for $\vecs = (s_1,s_2,0) \in \pa \R_+^3$.

\subsection{Expansion for small $a$}
We go one step further in the asymptotics of ${\cal T}$, by considering the regime of small radius $a$. 
The expression of  ${\cal T}$  involves  the maps
\begin{equation}
\label{}
\begin{array}{ccc}
\mathcal{T}_{i,j} : H^{-1/2}(\partial B) & \rightarrow &H^{1/2}(\partial B)\\
\vecf_j &\mapsto & \displaystyle \int_{\partial B} \vecK(\vecx_i+ a\cdot,\vecx_j + a\vecs) \, \vecf_j (\vecs) \,\mbox{d}\vecs\,,
\end{array}
\end{equation}
with  the Green kernel $\vecK$ given by Proposition \ref{prop:Green_function}:
$$
\vecK(\vecr,\vecr'):=\vecG(\vecr-\vecr')+\vecK_1(\vecr,\vecr')+\vecK_2(\vecr,\vecr')+\vecK_3(\vecr,\vecr') + \vecK_4(\vecr,\vecr').
$$
We recall that $\vecK_1, \vecK_2$ and $\vecK_3$ are defined in \eqref{DefK}, whereas  $\vecK_4$ is defined by (see \eqref{K1}):
$$
 \vecK_4(\vecr,\vecr') := - \eps \int_{\partial \R_+^3}\,  h(\vecs) \, \frac{\partial}{\partial z}\left(\vecs \mapsto \vecK^0(\vecs,\vecr) \right) \, \frac{\partial}{\partial z}\left(\vecs \mapsto \vecK^0(\vecs,\vecr') \right) \mbox{d}\vecs\,.
$$
%
%In order to proceed, we make an expansion of the vector fields of the dynamics and their Lie brackets with respect to $a$ (the radius of the balls) near 0.
%This part is devoted to the proof of the expansion of the Neumann to Dirichlet map (\ref{expansion}) together with its inverse (\ref{inverse}) at large arms' lengths. Let us first define for all $(i,j)\in \{1,2,3\}^2$, the linear map $\mathcal{T}_{i,j}$ as
%
%\begin{eqnarray*}
%\mathcal{T}_{i,j} : H^{-1/2}(\partial B) & \rightarrow &H^{1/2}(\partial B)\\
%\vecf_j &\mapsto & \int_{\partial B} \vecK(\vecx_i+ a\cdot,\vecx_j + a\vecs) \, \vecf_j (\vecs) \,\mbox{d}\vecs\,.
%\end{eqnarray*}
%We recall that the Green kernel $\vecK$ writes (following (\ref{DefK})) as
%$$
%\vecK(\vecr,\vecr')=\vecG(\vecr-\vecr')+\vecK_1(\vecr,\vecr')+\vecK_2(\vecr,\vecr')+\vecK_3(\vecr,\vecr')\,,
%$$
%where $\vecG$ is the Stokeslet (see (\ref{stokeslet})) and each kernel is given by the corresponding counterpart in (\ref{DefK}). 
Eventually, we call $\mathcal{T}^G$ the Neumann to Dirichlet map associated to $\vecG$
\begin{eqnarray*}
\mathcal{T}^G : H^{-1/2}(\partial B) & \rightarrow &H^{1/2}(\partial B)\\
\vecf &\mapsto & \int_{\partial B} \vecG( a(\cdot - \vecs)) \, \vecf (\vecs) \,\mbox{d}\vecs\,.
\end{eqnarray*}

\begin{proposition}
\label{propexp}
Let $(i,j)\in \{1,\cdots,N\}^2$. We have the following expansions, valid for $a \ll 1$:
\begin{itemize}
\item if $i\ne j$ then
\begin{equation}
\mathcal{T}_{i,j} = \vecK(\vecx_i,\vecx_j)   \langle \, \cdot \, , \vecI_d \rangle_{\partial B} +\vecR_1
\label{idifj}
\end{equation}
where $\ds ||\vecR_1||_{\mathcal{L}(H^{-1/2},H^{1/2})}=O\left(a\right)$\,,
\item otherwise
\begin{equation} 
\mathcal{T}_{i,i} = \mathcal{T}^G+\sum_{k=1}^4 \vecK_k(\vecx_i,\vecx_i)   \langle \, \cdot \, , \vecI_d \rangle_{\partial B} +\vecR_2
\label{ieqj}
\end{equation}
where $\ds ||\vecR_2||_{\mathcal{L}(H^{-1/2},H^{1/2})}=O\left(a\right)$\,.
\end{itemize}
\end{proposition}

\Proof
Let $(i,j)\in \{1,\cdots,N\}^2$ be such that $i\ne j$. For all  $\vecf_j \in H^{-1/2}(\partial B)$, we write
\begin{equation}
\left(\mathcal{T}_{i,j} - \vecK(\vecx_i,\vecx_j)\langle \cdot, \vecI_d \rangle\right)(\vecf_j)(\vecr) = \int_{\partial B}  \left(\vecK(\vecx_i+a\vecr,\vecx_j+a\vecs)-\vecK(\vecx_i,\vecx_j)\right)\vecf_j(\vecs)\mbox{d}\vecs\,.
\end{equation}
The point is that, as $i \neq j$, the kernel $\vecK$ is smooth in a neighborhood of $B_i \times B_j$. Hence, 
\begin{equation}
\left|\vecK(\vecx_i+a\vecr,\vecx_j+a\vecs) - \vecK(\vecx_i,\vecx_j)\right| = O\left(a\right)\,, \left| \na \vecK(\vecx_i+a\vecr,\vecx_j+a\vecs) - \vecK(\vecx_i,\vecx_j)\right| = O\left(a\right)
\end{equation}
 uniformly for $\vecr,\vecs \in B$. Estimate \eqref{idifj} follows straightforwardly. \\
%\begin{eqnarray*}
%&&\left|\nabla_{\vecr} \vecK(\vecx_i+a\vecr,\vecx_j+a\vecs) \right| = O\left(a\right)\,,\\
%&&\left|\nabla_{\vecs} \vecK(\vecx_i+a\vecr,\vecx_j+a\vecs) \right| = O\left(a\right)\,,\\
%&&\left|\nabla_{\vecr}\nabla_{\vecs} \vecK(\vecx_i+a\vecr,\vecx_j+a\vecs) \right| = O\left(a^2\right)\,.
%\end{eqnarray*}
%Therefore, we obtain $\forall \vecr \in B$
%\begin{eqnarray*}
%|\vecv_i(\vecr)|&\leq& \|\vecK(\vecx_i+a\vecr,\vecx_j+a\cdot)-\vecK(\vecx_i,\vecx_j)\|_{H^\frac12}\|\vecf_j\|_{H^{-\frac12}}\\
%&\leq& O\left(a\right)\|\vecf_j\|_{H^{-\frac12}}\,,
%\end{eqnarray*}
%and similarly
%\begin{eqnarray*}
%|\nabla_{\vecr}\vecv_i(\vecr)|&\leq& \|\nabla_{\vecr} \left(\vecK(\vecx_i+a\vecr,\vecx_j+a\cdot)\right)\|_{H^\frac12}\|\vecf_j\|_{H^{-\frac12}}\\
%&\leq& O\left(a\right)\|\vecf_j\|_{H^{-\frac12}}\,.
%\end{eqnarray*}
%This enables us to estimate the $H^\frac12$ norm of $\vecv_i$ on $\partial B$
%\begin{eqnarray*}
%\left\|\vecv_i\right\|_{H^\frac12(B)} &\leq& \left\|\vecv_i\right\|_{H^1(B)}\\
% &=& \left(\left\|\vecv_i\right\|^2_{L^2(B)}+\left\|\nabla \vecv_i\right\|^2_{L^2(B)}\right)^\frac12\\
% &\leq&O\left(a\right)\|\vecf_j\|_{H^{-\frac12}}\,,
%\end{eqnarray*}
%which proves (\ref{idifj}).

The proof of  (\ref{ieqj}) is similar:  we have for all 
$\vecf_i \in H^{-1/2}(\partial B)$
\begin{equation}
\left(\mathcal{T}_{i,i} - \mathcal{T}^G -   \vecK(\vecx_i,\vecx_j)\langle \cdot, \vecI_d \rangle\right)(\vecf_i)(\vecr) = \int_{\partial B}  \sum_{k=1}^4\left(\vecK_k(\vecx_i+a\vecr,\vecx_i+a\vecs)-\vecK_k(\vecx_i,\vecx_i)\right)\vecf_i(\vecs)\mbox{d}\vecs\,,
\end{equation}
where none of the $\vecK_k$'s is singular near $B_i \times B_i$. 
%
%use the decomposition (\ref{DefK}) where none of the kernels $(\vecK_i)_{i=1,2,3}$ is singular. Therefore $\forall \vecr\in \partial B$
%\begin{eqnarray*}
%\vecu_i(\vecr):=( \mathcal{T}_{i,i} \vecf_i)(\vecr) &=& \int_{\partial B} \vecK(\vecx_i+a\vecr,\vecx_i+a\vecs)\vecf_i(\vecs)\mbox{d}\vecs\,\\
%&=& \int_{\partial B} \vecG(a(\vecr-\vecs))\vecf_i(\vecs)\mbox{d}\vecs+\int_{\partial B} (\vecK_1+\vecK_2+\vecK_3+\vecK_4)(\vecx_i+a\vecr,\vecx_i+a\vecs)\vecf_i(\vecs)\mbox{d}\vecs\\
%&=&\mathcal{T}^G \vecf_i + \int_{\partial B} (\vecK_1+\vecK_2+\vecK_3 + \vecK_4)(\vecx_i+a\vecr,\vecx_i+a\vecs)\vecf_i(\vecs)\mbox{d}\vecs\,.
%\end{eqnarray*}
%We finish as before, having remarked that for $l=1,\dots,4$
%\begin{equation}
%\vecK_l(\vecx_i+a\vecr,\vecx_i+a\vecs) = \vecK_l(\vecx_i,\vecx_i) + O\left(a\right)\,.
%\end{equation}
\EndProof
%\begin{figure}[ht]
%\centering
%\includegraphics[scale=0.4]{Figures/Regime}
%\caption{\label{Regime}}
%\end{figure}
As a simple consequence of the previous propositions, we have
\begin{proposition}
For every $\vecf \in H^{-1/2}_N$, for all $(\vecx,\vecxi)$, 
\begin{equation}
\left(\mathcal{T} \vecf \right)_i (\vecr) = \mathcal{T}^G \vecf_i + \sum_{l=1}^4 \vecK_l(\vecx_i,\vecx_i)   \langle \vecf_i, \vecId \rangle_{\partial B}  + \displaystyle \sum_{j\neq i} \vecK(\vecx_i,\vecx_j) \langle \vecf_j, \vecId\rangle_{\partial B} + \mathcal{R}_i(\vecf),
\label{expansion}
\end{equation}
with $\ds \Vert \mathcal{R}_i\Vert_{\mathcal{L}(H_N^{-1/2},H_N^{1/2})} =O\left(a + \eps^2 \right)$, and $i=1...N$.
\end{proposition}
\Proof
By Proposition  \ref{prop:Green_function}: for all $i=1...N$, and all $\vecr\in \partial B$
\begin{eqnarray*}
\left(\mathcal{T} \vecf \right)_i (\vecr)&:=& \displaystyle \int_{\partial B} \vecK(\vecx_i +a\vecr, \vecx_i + a\vecs) \, \vecf_i(\vecs) \mbox{ds} + \sum_{i \neq j} \int_{\partial B} \vecK(\vecx_i +a\vecr, \vecx_j + a\vecs) \, \vecf_j(\vecs) \mbox{d}\vecs  + {\cal R}^\eps(\vecf) \\
&=& \mathcal{T}_{i,i}\vecf_i + \sum_{j \neq i} \mathcal{T}_{i,j} \vecf_j + {\cal R}^\eps(\vecf), \qquad \| {\cal R}^\eps \|_{\mathcal{L}(H_N^{-1/2},H_N^{1/2})} = O(\eps^2) 
%\label{Expressionu_i}
\end{eqnarray*}
and the result follows from the application of (\ref{idifj}) and (\ref{ieqj}) of Proposition \ref{propexp}.
\EndProof
\begin{proposition}
For every $\vecu \in H_N^{1/2}$, for all $(\vecp,\vecxi)$, one has 
\begin{equation}
\label{inverse}
\begin{array}{ll}
\left(\mathcal{T}^{-1} \vecu \right)_i = & ({\mathcal{T}^G})^{-1} \left( \vecu_i - \displaystyle\sum_{k=1}^4  \vecK_k(\vecx_i,\vecx_i) \langle ({\mathcal{T}^G})^{-1} \vecu_i, \vecId \rangle_{\partial B} \right) -\\
& \hspace*{2cm} \displaystyle ({\mathcal{T}^G})^{-1} \left(  \sum_{j\neq i} \vecK(\vecx_i,\vecx_j) \langle ({\mathcal{T}^G})^{-1} \vecu_j, \vecId \rangle_{\partial B}\right) + \tilde{\mathcal{R}}_i(\vecu)
\end{array}
\end{equation}
with $\ds \Vert \tilde{\mathcal{R}}_i\Vert_{\mathcal{L}(H_N^{1/2},H_N^{-1/2})} =O\left(a^3 + a^2 \eps^2 \right)$, $i=1...N$.
\end{proposition}

\Proof
We recall that
\begin{equation*}
\mathcal{T}^G : H^{-\frac12}(\partial B) \rightarrow H^{\frac12}(\partial B), \quad 
\vecf \mapsto  \int_{\partial B} \vecG(a(\cdot-\vecs)) \vecf(\vecs)\,d\vecs\,,
\end{equation*}
and define for $l=1,\dots,4$ the operators
\begin{equation*}
\mathcal{S}_l : H^{-\frac12}(\partial B) \rightarrow H^{\frac12}(\partial B), \quad 
\vecf \mapsto \int_{\partial B} \vecK_l(\vecx_i,\vecx_i) \vecf(\vecs)\,d\vecs\,,
\end{equation*}
and eventually
\begin{equation*}
\mathcal{S}_{i,j} : H^{-\frac12}(\partial B) \rightarrow H^{\frac12}(\partial B), \quad 
\vecf \mapsto \int_{\partial B} \vecK(\vecx_i,\vecx_j) \vecf(\vecs)\,d\vecs\,.
\end{equation*}
Notice that for all $\vecf\in H^{-\frac12}(\partial B)$, $\mathcal{S}_l \vecf$ and $\mathcal{S}_{i,j}\vecf$ are  constant applications. \\

That these operators are continuous operators from $H^{-\frac12}(\partial B)$ into $H^{\frac12}(\partial B)$ is classical. We  are only interested in estimating their norms, and more precisely in  the way they depend on $a$ in the limit $a\rightarrow 0$. Notice that since the kernel $\vecG$ is homogeneous of degree -1, one has
\begin{equation}
\|\mathcal{T}^G\|_{\mathcal{L}(H^{-1/2},H^{1/2})} =O\left(\frac{1}{a}\right)\,\mbox{ and }\left\|\left(\mathcal{T}^G\right)^{-1}\right\|_{\mathcal{L}(H^{1/2},H^{-1/2})} =O\left(a\right)\,.
\end{equation}
As far as $\mathcal{S}_l$ is concerned, we get that (since $\ds |\vecK_l(\vecx_i,\vecx_i)| = O\left(1\right)$)
\begin{equation}
\|\mathcal{S}_l\|_{\mathcal{L}(H^{-1/2},H^{1/2})} =O\left(1\right)\,,
\end{equation}
and similarly
\begin{equation}
\|\mathcal{S}_{i,j}\|_{\mathcal{L}(H^{-1/2},H^{1/2})} =O\left(1\right)\,.
\end{equation}
When $a \rightarrow 0$ this enables us to invert (\ref{expansion}) leading to (\ref{inverse}).
\EndProof

\section{Controllability of the Three-sphere swimmer}
\label{sec:Three_sphere}

%\subsection{Asymptotic expansion of the motion equation}
 %\label{Subsec:Self-propulsion}

We deal in this section  with the controllability of the 3-sphere swimmer, namely Theorem \ref{thm:three_sphere}. 
%This theorem is established in several steps. First, we describe the particular dynamics associated with the 3-sphere swimmer. Then, we give an asymptotic expansion of the vector fields of the motion equation. We  use this asymptotic expansion to determine the dimension of the Lie algebra "almost everywhere" (i.e., for almost every radius of the ball, position of the swimmer and amplitude of the roughness). This allows us to conclude with Chow's theorem.
\\

\subsection{Preliminary remarks on the 3-sphere dynamics}

%In what follows, the swimmer is described by the coordinates for the middle sphere $\vecX:=(x,y,z)$ and the angles $(\theta, \phi)$ represented in the Fig. \ref{Fig:coordinates}. 

%In what follows, we show that the rotation of the swimmer around its axis appears at order too small to be considered in the expansion of its motion equation. The parameters used to describe the motion of the swimmer is the coordinate of the middle sphere, noted by $(x,y,z)$, and the angles $(\theta, \phi)$ represented in the Fig. \ref{Fig:coordinates}.
We must first come back to equation \eqref{eq:dynamics1} \eqref{eq:dynamics}, in the particular case of the 3-sphere swimmer.  Remember that the writing in this equation was slightly abusive: we had denoted by $\dot{\vecp}$ the vector $\left( \begin{smallmatrix} \Om \\ \vecv \end{smallmatrix} \right)$ associated to the rigid movement of the swimmer, see \eqref{rigidvf}. In our case, $\Omega = \left( \begin{smallmatrix} \Omega_1 \\ \Omega_2 \\ \Omega_3 \end{smallmatrix} \right)$ and $\vecv = \dot{\vecx}_c =  \left( \begin{smallmatrix} v_1 \\ v_2 \\ v_3 \end{smallmatrix} \right)$ 
 are  respectively the angular  velocity and the linear velocity of the middle sphere, decomposed in an arbitrary orthonormal basis $(\vece_i)$. Moreover, it is natural to take for $\vece_1$ the unit vector of the 3-sphere axis. Let $\theta$ be the angle between the swimmer's axis and $\vece_z$, while $\phi$ is the angle between the $x$-axis and the projection of the swimmer in $Oxy$ plane (see figure \ref{Fig:coordinates_intro}). Then, the unit vector of the 3-sphere axis reads (in the canonical basis)
 $\vece_1=\left(\begin{array}{c} \cos(\phi)\sin(\theta) \\ \sin(\phi)\sin(\theta) \\ \cos(\theta) \end{array}\right)$. It is completed into an orthonormal basis by defining 
 \begin{displaymath}
\vece_2=\left(\begin{array}{c} \cos(\phi)\cos(\theta) \\ \sin(\phi)\cos(\theta) \\ -\sin(\theta) \end{array}\right)\,, \quad \quad \quad \vece_3 = \left(\begin{array}{c} -\sin(\phi) \\ \cos(\phi) \\ 0 \end{array}\right)\,.
\end{displaymath}
Hence, a rigorous writing of \eqref{eq:dynamics1} or \eqref{eq:dynamics} is 
\begin{equation} \label{rigorousdyn}
\vecM \left( \begin{smallmatrix} \Om \\ \vecv \end{smallmatrix} \right) + \vecN = 0, \quad \textrm{or} \: \left( \begin{smallmatrix} \Om \\ \vecv \end{smallmatrix} \right) = -\vecM^{-1} \vecN. 
\end{equation}
A crucial remark is that $\vecM$ and $\vecN$ do not depend on the whole of $\vecp$. The angle $\theta_1$ of rotation around the swimmer's axis is not involved, as it is irrelevant to the swimmer's position, orientation or elongation. In particular, keeping only the five bottom lines of the last system, we end up with a closed relation of the type
\begin{equation}  \label{reducedode}
\left( \begin{smallmatrix} \dot\theta_2 \\ \dot\theta_3 \\ \dot\vecx_c \end{smallmatrix} \right) = \sum_{i=1}^2 \tilde \vecF_i\left( \left( \begin{smallmatrix} \theta_2 \\\theta_3 \\ \vecx_c \end{smallmatrix} \right) \right) \dot{\xi_i} 
\end{equation}
where $\theta_2$ and $\theta_3$ are the rotation angles around $\vece_2$ and $\vece_3$ respectively. Then, by the analyticity of the $\tilde \vecF_i$'s and Chow's theorem, it remains to prove that there exists some $(\eps,a,\theta_2,\theta_3,\vecx_c)$ such that 
$$ \dim Lie_{(\theta_2,\theta_3,\vecx_c)}\left( \left( \begin{smallmatrix} 1 \\ 0 \\ \tilde \vecF_1 \end{smallmatrix} \right), \; \left( \begin{smallmatrix} 0 \\ 1 \\ \tilde \vecF_2 \end{smallmatrix} \right) \right) \: = \: 7. $$
Actually, we shall not work directly with angles $\theta_2, \theta_3$. We find it more convenient to work with the angles $\theta, \phi$ introduced above (see Figure \ref{Fig:coordinates_intro}). From the relation 
$\displaystyle \frac{d}{dt} \vece_1 = \Omega \times \vece_1$, we infer that 
$$ \Omega_2 = -\sin\theta \dot{\phi}, \quad \Omega_3 = \dot{\theta}.  $$ 
Note that in the special case $\sin \theta = 0$, the angle $\phi$ coincides with the useless angle $\theta_1$. Moreover, the mapping $(\theta_2,\theta_3) \rightarrow (\theta, \phi)$ is not a diffeomorphism in the vicinity of $\theta \equiv 0 [\pi]$. {\em Thus, we shall restrict to orientations of the swimmer for which 
\begin{equation} \label{sinneq0}
|\sin \theta | \ge \delta > 0.
\end{equation}}
We shall establish the maximality of the Lie algebra at points satisfying this condition. 
\\

Before entering the computation of this Lie algebra, we state a technical lemma, that will somehow allow us to neglect the rotation around the swimmer's axis. As mentioned before, we assume inequality \eqref{sinneq0}. We have 
\begin{lemma}
\label{lemma:estimationOmega}
There exists a constant $C$ which does not depend on $a$ and $\epsilon$ such that
$$
|\Omega_1| \leq C \, \left( |\dot{\theta}| + |\dot{\phi}| + |\dot{\vecx}_c| +  |\dot{\vecxi}| \right)\,.
$$
\end{lemma}

\Proof
We go back to the first identity in \eqref{rigorousdyn}. The first line gives
\begin{equation}
\label{eq:expressionOmega}
\vecM_{1,1} \,\Omega_1 =  - \vecN_1 + \vecM_{1,2}\, \sin(\theta) \dot{\phi} - \vecM_{1,3}\, \dot{\theta} - \vecM_{1,4}\, v_1  - \vecM_{1,5}\, v_2 - \vecM_{1,6}\, v_3 \,.
\end{equation}
We recall that, in the definitions of $\vecM$ and $\vecN$, we denoted by $\vecu_i$ and $\vecu^d$ some solutions of the Stokes equation, with zero Dirichlet condition at the wall, and inhomogeneous Dirichlet conditions at the ball.  The Dirichlet data is $\vece_i \times (\vecx - \vecx_c)$ for $i=1,2,3$,  $\vece_{i-3}$ for $i=4,5,6$, and $\vecu_d$ for $\vecu^d$. In the case of the 3-sphere swimmer,  $\vecu_d$ is 
$-\dot{\xi}_1 \vece_1$ on the sphere $\pa B_1$, $0$  on the middle sphere and $\dot{\xi}_2 \vece_1$ on the sphere $\pa B_3$. 
\\

Let us first examine 
\begin{equation}
\label{eq:expressionM}
\begin{aligned}
\vecM_{1,1} \: & = \:  
\sum_{l=1}^3 \int_{\partial B} \left(\vecx_l - \vecx_c + a\vecr \right) \times \vece_1  \cdot \left( \mathcal{T}^{-1}( \vece_1 \times a \vecr, \vece_1 \times a \vecr, \vece_1 \times a \vecr) \right)_l    \mbox{d} \sigma  \\
& = \: 3 \int_{\pa B} a \vecr \times \vece_1 \cdot  \left(\mathcal{T}^{-1}( \vece_1 \times a \vecr, \vece_1 \times a \vecr, \vece_1 \times a \vecr)  \right)_l       \mbox{d} \sigma  
\end{aligned}
\end{equation}
using that $(\vecx_l - \vecx_c) \times \vece_1 = 0$. We then use the expansion  \eqref{inverse}. We recall the well-known fact that the rotation are eigenfunctions of  $\left(\mathcal{T}^G\right)^{-1}$, with associated eigenvalue $3\mu a$. In particular,
$$  
\left(\mathcal{T}^G\right)^{-1}(\vece_1 \times a\vecr) \: = \:  3 \mu a \vece_1 \times a\vecr,  \quad \textrm{and} \: \langle \left(\mathcal{T}^G\right)^{-1}(\vece_1 \times a\vecr) , \vecId \rangle_{\partial B} = 0.$$
We find then  easily that 
$ \vecM_{1,1} =  -3 \mu a^3 \: + \: O(a^5 + \eps^2 a^4) $. \\

Then, we examine  
\begin{equation*}
\vecN_{1}   =   \sum_{l=1}^3  \int_{\partial B} \left(\vecx_l - \vecx_c + a\vecr \right) \times \vece_1 \cdot \left( \mathcal{T}^{-1}(-\dot\xi_1 \, \vece_1, 0, \dot\xi_2 \, \vece_2) \right)_l \,\mbox{d}\sigma.
\end{equation*}
Again, we can expand $\mathcal{T}^{-1}$ using \eqref{inverse}. This time, we use that translations are  eigenfunctions of $\left(\mathcal{T}^G\right)^{-1}$ with associated eigenvalue $\frac{3}{2}\mu a$. Thus, 
$$ \left(\mathcal{T}^G\right)^{-1}(\vece_1) \: = \:  \frac{3}{2}\mu a \,  \vece_1. $$
It follows that the first terms in the expansion vanish, and  we find 
$$ \vecN_1 = O((a^4 + a^3 \eps^2) \,  |\dot\vecxi|) $$
\\

The remaining terms $M_{1,j}$, $j=2,...,4$ can be handled with similar arguments. The  lemma follows straightforwardly.

\EndProof

\subsection{Asymptotics of the 3-sphere dynamics}
 \label{Subsec:Self-propulsion}
We shall now provide an accurate description of the 3-sphere dynamics: broadly speaking, the point is to obtain an explicit expansion of the $\tilde{\vecF}_i$'s in \eqref{reducedode} (with angles $\theta_2,\theta_3$ replaced by $\theta,\phi$, see remark above). We remind that the dynamics (that is the 6x6 system in \eqref{rigorousdyn})  is governed by  self-propulsion: it corresponds to 
\begin{itemize}
\item The sum of the forces on the swimmer being zero. 
\item The sum of the torques on the swimmer being zero. 
\end{itemize}

{\em Forces.} By the definition of the swimmer, each sphere obeys a rigid body motion. More precisely, the velocity of each point $\vecr$ of the $l$th sphere expresses as a sum of a translation and a rotation as
\begin{equation}
\vecu^S_l(\vecr)=\vecu_{T_l} + \vecu_{R_l}(\vecr)\,,
\label{mvtsolide}
\end{equation}
where $\vecu_{T_l}$ is constant on $\partial B$ while $\vecu_{R_l}(\vecr)=\Omega \times a\vecr$ (remember that all quantities are expressed on the unit sphere $\partial B$). The vanishing of the total force, due to self-propulsion, reads
\begin{equation}
\label{SommeDesForces}
\begin{array}{ll}
\displaystyle \sum_{l} \int_{\partial B} \vecf_l = \displaystyle \sum_{l} \int_{\partial B} \left( \mathcal{T}^{-1}\left( \vecu^S_1, \vecu^S_2, \vecu^S_3 \right) \right)_l = 0\,.
\end{array}
\end{equation}
Plugging (\ref{mvtsolide}) in (\ref{SommeDesForces}) and using (\ref{inverse}) leads to
\begin{multline}
\label{TermesDiagonaux}
 \sum_l \int_{\partial B}
({{\mathcal{T}}^G})^{-1} \left({\vecu_T}_l + {\vecu_R}_l - \sum_{k=1}^4  \vecK_k(\vecx_l,\vecx_l) \langle ({{\mathcal{T}}^G})^{-1} ({\vecu_T}_l + {\vecu_R}_l), \mbox{Id} \rangle_{\partial B} \right) -  \\ 
 ({{\mathcal{T}}^G})^{-1} \left(  \sum_{j\neq i} \vecK(\vecx_i,\vecx_j) \langle (\mathcal{T}^G)^{-1} ({\vecu_T}_l + {\vecu_R}_l), \vecId \rangle_{\partial B}\right) = \left(O\left(a^3\right)+ O\left(a^2\epsilon^2\right)\right)||\vecu||\,.
\end{multline}
where $\| \vecu \| = \| (\vecu_i^S) \|$ is any norm on the n-uplets of rigid vector fields over the ball. Here, 
\begin{equation}
\| u \| \:  = \: O(| \dot\theta| + |\dot\phi| + | \dot \vecx_c |  + | \Omega_1 | )  \: = \: O(| \dot\theta| + |\dot\phi| + | \dot \vecx_c |) 
\label{conseqlemma}
\end{equation}
where the last equality comes from Lemma \ref{lemma:estimationOmega}. 	   
As mentioned earlier, it is well known that both translations and rotations are eigenfunctions of the Dirichlet to Neumann map of the three dimensional Stokes operator outside a sphere. Namely
$$
\left({{\mathcal{T}}^G}\right)^{-1} {\vecu_T}_l = \lambda_T{\vecu_T}_l\mbox{ and }  \left({{\mathcal{T}}^G}\right)^{-1}{\vecu_R}_l = \lambda_R {\vecu_R}_l\,.
$$
It is also well-known that $\lambda_T=\frac{3\mu a}{2}$, $\lambda_R =  3\mu a$, leading in particular to the celebrated Stokes formula
$$
\int_{\partial B}
\left({{\mathcal{T}}^G}\right)^{-1} \vecu_{T_l} \mbox{d}\vecs= 6\pi\mu a \,\vecu_{T_l}
$$
 We also remark that due to
$
\int_{\partial B}  \vecu_{R_l} \mbox{d}\vecs = 0\,$, we have
$
\int_{\partial B}
\left({{\mathcal{T}}^G}\right)^{-1} \vecu_{R_l} \mbox{d}\vecs = 0\,.
$
We therefore obtain
\begin{multline}
\displaystyle 6\pi\mu a \sum_l \left({\vecu_T}_l - 6\pi\mu a\sum_{k=1}^4  \vecK_k(\vecx_l,\vecx_l)   {\vecu_T}_l  - 6\pi\mu a \sum_{j\neq i} \vecK(\vecx_l,\vecx_j) {\vecu_T}_j \right) \\= \left(O\left(a^3\right)+ O\left(a^2\epsilon^2\right)\right)||\vecu||\,.
\label{TotalForce}
\end{multline}

{\em Torques.}
We now compute the torque with respect to the center $\vecx_c$ of the middle ball $B_2$. Self-propulsion of the swimmer implies that this torque vanishes:
\begin{equation}
\label{TotalTorque}
0 =  \int_{\partial B} (\vecx_1-\vecx_2 + a\vecr) \times \vecf_1(\vecr)+\int_{\partial B} a\vecr \times \vecf_2(\vecr)+\int_{\partial B} (\vecx_3-\vecx_2 + a\vecr) \times \vecf_3(\vecr)
= \vecI_1 + \vecI_2 + \vecI_3\,,
\end{equation}
with the quantities $\vecI_1$, $\vecI_2$ and $\vecI_3$  given below.
\begin{eqnarray*}
\vecI_1 &=& \int_{\partial B} (\vecx_1-\vecx_2+a\vecr) \times \vecf_1(\vecr)=  \int_{\partial B} (-\xi_1\vece_1+a\vecr) \times\left( \mathcal{T}^{-1}\left( \vecu^S_1, \vecu^S_2, \vecu^S_3 \right) \right)_1
\\
&=& \int_{\partial B} \left( -\xi_1 \vece_1 + a\vecr \right) \times   ({{\mathcal{T}}^G})^{-1}\left({\vecu_T}_1 + {\vecu_R}_1 - 6\pi \mu a \sum_{k=1}^4  \vecK_k(\vecx_1,\vecx_1) {\vecu_T}_1\right.\\
&&\hspace*{5cm}\left. - 6\pi \mu a \sum_{j\neq 1} \vecK(\vecx_1,\vecx_j) {\vecu_T}_j +O\left(a^2+ a \eps^2\right ) ||\vecu||\right) \\
&=&-6\pi \mu a  \xi_1 \vece_1 \times   \left({\vecu_T}_1 - 6\pi \mu a \sum_{k=1}^4  \vecK_k(\vecx_1,\vecx_1) {\vecu_T}_1 - 6\pi \mu a \sum_{j\neq 1} \vecK(\vecx_1,\vecx_j) {\vecu_T}_j\right)\\
& &\hspace*{5cm}+ \left(O\left(a^3\right)+O\left(a^2\epsilon^2\right)\right) ||\vecu||\,.\\
\end{eqnarray*}
Similarly, we get,
\begin{eqnarray*}
\vecI_2 &=& a\int_{\partial B} \vecr \times \vecf_2(\vecr)=  a\int_{\partial B} \vecr \times\left( \mathcal{T}^{-1}\left( \vecu^S_1, \vecu^S_2, \vecu^S_3 \right) \right)_2\\
&=& a\int_{\partial B} \vecr \times   ({{\mathcal{T}}^G})^{-1} \left({\vecu_T}_2 + {\vecu_R}_2 - 6\pi \mu a \sum_{k=1}^4  \vecK_k(\vecx_2,\vecx_2) {\vecu_T}_2 \right.\\
&&\hspace*{5cm}\left.- 6\pi \mu a \sum_{j\neq 2} \vecK(\vecx_2,\vecx_j) {\vecu_T}_j +O\left(a^2 + a \eps^2\right)||\vecu||\right) \\
&=& \:  O\left(a^3 + a^3\epsilon^2\right) ||\vecu||\,.
\end{eqnarray*}
Finally,
\begin{eqnarray*}
\vecI_3 &=& \int_{\partial B} (\vecx_3-\vecx_2+a\vecr) \times \vecf_3(\vecr)\\
&=&6\pi \mu a  \xi_2 \vece_1 \times   \left({\vecu_T}_3- 6\pi \mu a \sum_{k=1}^4  \vecK_k(\vecx_3,\vecx_3) {\vecu_T}_3  - 6\pi \mu a \sum_{j\neq 3} \vecK(\vecx_3,\vecx_j) {\vecu_T}_j\right)\\
&&\hspace*{5cm}+ \left(O\left(a^3\right)+O\left(a^2\epsilon^2\right)\right)||\vecu||\,.\\
\end{eqnarray*}
Denoting by $\vecA$ the matrix
\begin{equation}
\vecA=\left(
\begin{array}{ccc}
\vecA_{11} & \vecA_{12}  & \vecA_{13}\\
\vecA_{21} & \vecA_{22}  & \vecA_{23}\\
\vecA_{31} & \vecA_{32}  & \vecA_{33}
\end{array}
\right)
\label{defA1}
\end{equation}
where for $i=1,2,3$
\begin{equation}
\vecA_{ii}=\vecId - 6\pi\mu a \sum_{l=1}^4 \vecK_l(\vecx_i,\vecx_i)
\label{defA2}
\end{equation}
and for $i,j=1,2,3$ with $i\neq j$
\begin{equation}
\vecA_{ij}=-6\pi\mu a \,\vecK(\vecx_i,\vecx_j)
\label{defA3}
\end{equation}
and $\vecS$ the matrix
$$
\vecS
=
\left(
\begin{array}{ccc}
\vecId & \vecId & \vecId \\
 -\xi_1 \vece_1 \times  & 0 & +\xi_2 \vece_1 \times\\
\end{array}
\right)\,,
$$
we can rewrite the self-propulsion assumption (\ref{TotalForce}), (\ref{TotalTorque}) as 
\begin{equation}
\vecS \vecA \left(
\begin{array}{c} \vecu_{T_1}\\ \vecu_{T_2} \\  \vecu_{T_3}
\end{array}\right) =\left(O\left(a^2\right)+O\left(a\epsilon^2\right)\right)||\vecu||.
\label{TotalForceTorque}
\end{equation}
Terms involving the $\vecu_{R_l}$'s are included in the r.h.s. 
\\

We now express $\vecu_{T_1}, \vecu_{T_2}$ and $\vecu_{T_{3}}$ in terms of
$\dot\vecx_c,\dot{\theta},\dot{\phi}$ and $\dot{\vecxi}$. 
Since $\vecu_{T_2}$ is the velocity of the center of the ball $B_2$, one has

$$
\vecu_{T_2} = \dot \vecx_c = \left(\begin{array}{c} \dot{x} \\ \dot{y} \\ \dot{z} \end{array}\right)\quad \textrm{in the canonical basis of} \: \R^3.
$$
Then, by using $
\frac{d}{dt}\vece_1= \dot{\theta}\,\vece_2+\sin(\theta)\dot{\phi}\,\vece_3\,,
$ we get
$$
\vecu_{T_1} = \vecu_{T_2} - \xi_1\left(\dot{\theta} \vece_2 + \sin(\theta)\dot{\phi}\,\vece_3 \right) - \dot{\xi}_1 \vece_1\,, \quad \vecu_{T_3} = \vecu_{T_2} + \xi_2 \left(\dot{\theta} \vece_2 + \sin(\theta)\dot{\phi}\,\vece_3\right) + \dot{\xi}_2 \vece_1\,\,.
$$
In matrix form, all this reads 
Then, the speed $\vecu_{T_i}$ ($i=1,2,3$) is expressed as
\begin{equation}
\left(\begin{array}{c} \vecu_{T_1}\\ \vecu_{T_2} \\  \vecu_{T_{3}}
\end{array}\right) = \vecT \,  \left(\begin{array}{c}\Omega_1\\ \dot{\theta} \\ \dot{\phi}  \\ \dot{x} \\ \dot{y} \\ \dot{z} \end{array}\right)+\vecU \,\dot{\vecxi} \,.
\label{uti2}
\end{equation}

\noindent with
$$
\vecT= \left(
\begin{array}{cccc}
0 & - \xi_1 \vece_2 & - \xi_1 \sin(\theta)\vece_3 & \vecId \\
\vdots& \begin{array}{c} 0 \\ 0 \\0 \end{array} &  \begin{array}{c} 0 \\ 0 \\0 \end{array} &\vecId  \\
0 & +\xi_{2} \vece_2 & +\xi_{2}\sin(\theta)\vece_3&\vecId 
\end{array}
\right)\,, \quad \textrm{and} \quad 
\vecU= \left(
\begin{array}{cc}
-\vece_1 & 0 \\
0 & \vdots \\
\vdots &  0\\
0& \vece_1 \\
\end{array}
\right).
$$
Combining with \eqref{TotalForceTorque}, the motion equation \eqref{rigorousdyn} becomes
\begin{equation}
 \left(\vecS \vecA+\vecR_1\right)\left( \vecT \, \left(\begin{array}{c}\Omega_1\\ \dot{\theta} \\ \dot{\phi}  \\ \dot{x} \\ \dot{y} \\ \dot{z} \end{array}\right)+ (\vecU + \vecR_2) \,\dot{\vecxi}\right)
 = 0
\label{EquationMouvement2}
 \end{equation}
where the residual matrices $\vecR_1, \vecR_2$ satisfy 
$$
|\vecR_1|  + |\vecR_2| = \left(O\left(a^2\right)+O(a\epsilon^2)\right)\,
$$
using \eqref{conseqlemma}. 
Finally, we only keep the five bottom lines of this system. It yields the following 5x5 system
\begin{equation}
 \left(\tilde{\vecS} \vecA+\tilde{\vecR}\right)\left( \tilde{\vecT} \, \left(\begin{array}{c} \dot{\theta} \\ \dot{\phi}  \\ \dot{x} \\ \dot{y} \\ \dot{z} \end{array}\right)+ \vecU \,\dot{\vecxi}\right)
 = 0\,,
\label{eq:AsymptoticDynamicsThreeSphere}
 \end{equation}
where 
$$
\tilde{\vecS} \: := \: \left( S_{i,j} \right)_{2 \le i \le 6, 1 \le j \le 9}, \quad  \tilde{\vecT} \: := \:  \left(
\begin{array}{ccc}
 - \xi_1 \vece_{2} & - \xi_1 \sin(\theta)\vece_3 & \vecId \\
 \begin{array}{c} 0 \\ 0 \\0 \end{array} &  \begin{array}{c} 0 \\ 0 \\0 \end{array} &\vecId  \\
 +\xi_{2} \vece_2 & +\xi_{2}\sin(\theta)\vece_3&\vecId 
\end{array}
\right)\,, 
$$
and where  the residual matrices still satisfy  $| \tilde{\vecR}_1 | + | \tilde \vecR_2|  = O\left(a^2\right)+O(a\epsilon^2)$. We leave to the reader to check that 
$\tilde \vecS \vecA \tilde \vecT = \tilde \vecS \tilde \vecT + O(a)$ is invertible, with $|(\tilde \vecS \vecA \tilde \vecT)^{-1}| = O(1)$ uniformly in $a$ and $\eps$. Then, we can write  system \eqref{eq:AsymptoticDynamicsThreeSphere} as
 \begin{equation} \label{dynamicsinvert}
  \left( \begin{smallmatrix} \dot\theta \\ \dot\phi \\ \dot x \\ \dot y \\ \dot z \end{smallmatrix} \right)  = - 
  (\tilde \vecS \vecA \tilde \vecT)^{-1} \tilde \vecS \vecA  \vecU \dot\vecxi \: + \: \tilde \vecR \dot \vecxi
 \end{equation}
with $| \tilde \vecR|  = O(a^2 + \eps^2 a)$. 

\subsection{Reachable set}
We are now ready to prove Theorem \ref{thm:three_sphere}. We drop the tilda in the ODE \eqref{dynamicsinvert} and express it as 
\begin{equation} \label{eq:dynamicsThreeSphereControleTools}
\dot{\vecX} \: = \:   \vecF_1(\vecX) \dot\xi_1 \: + \: \vecF_2(\vecX)  \dot\xi_2, \quad \vecX \: := \:  \left( \begin{smallmatrix} \xi_1 \\ \xi_2 \\ \theta \\ \phi \\x \\ y \\ z \end{smallmatrix} \right). 
\end{equation}
To expand the $\vecF_i$'s, we decompose the matrix $\vecA$  into three matrices:  
$\vecA\:  := \:  \vecId +  \vecA^1 + \vecA^2$
where 
\begin{equation*}
\vecA^1_{ii}= -6\pi\mu a \sum_{k=1}^3 \vecK_k(\vecx_i,\vecx_i)  \:\:  \forall \, i\quad 
\vecA^1_{ij}=-6\pi\mu a \,\left(\vecG(\vecx_i,\vecx_j) + \sum_{k=1}^3 \vecK_k(\vecx_i,\vecx_j)\right) \:\: \forall \, i \neq j
\end{equation*}
and where
\begin{equation*}
\vecA_{i,j}^2= -6\pi\mu a\, \vecK_4(\vecx_i,\vecx_j) \quad \forall \, i,j. 
\end{equation*}
%These three vector fields are expressed
Thanks to \eqref{eq:AsymptoticDynamicsThreeSphere}, we get an expansion of the form  $\vecF_i  \: := \vecF_i^0 + \vecF_i^1+ \vecF_i^2+ \vecR_i$ where $\vecF_i^0$, $\vecF^1_i$ and $\vecF_i^2$ are respectively  the zero order term, the term of order $	a$ and the term of order $\eps a$. The remainder is    $\vecR_i =   \left(O\left(a^2\right)+O(a\epsilon^2)\right)$.
%we get an expansion of the vector fields of the dynamics \eqref{eq:AsymptoticDynamicsThreeSphere}, $\vecF_i$ ($i=1,2$), as the sum $\vecF_i = \vecF_i^0 + \vecF_i^1+ \vecF_i^2+ \vecR_i$ where $\vecF_i^0$ is the term of order zero, $\vecF_i^1$ (resp.  $\vecF_i^2$) is estimated as $O\left(a\right)$ (resp. $O\left(\epsilon a\right)$) and $\vecR_i =   \left(O\left(a^2\right)+O(a\epsilon^2)\right)$.
%\\
These vector fields are given by
%the first order term derived from the roughness.
\begin{equation}
\label{eq:vectorfields}
\begin{aligned}
\vecF^0_i & = \left( \begin{smallmatrix} \vece_i \\    -(\vecS\vecT)^{-1} (\vecS\vecU)  \vece_i \end{smallmatrix} \right),  \\
\vecF^1_i & =  \left( \begin{smallmatrix} 0 \\  \left((\vecS\vecT)^{-1}\vecS\vecA^1\vecT(\vecS\vecT)^{-1}\vecS\vecU - (\vecS\vecT)^{-1} \vecS \vecA^1 \vecU \right)  \vece_i\end{smallmatrix} \right), \\
\vecF^2_i & = \left( \begin{smallmatrix} 0 \\  \left((\vecS\vecT)^{-1}\vecS\vecA^2\vecT(\vecS\vecT)^{-1} \vecS\vecU- (\vecS\vecT)^{-1} \vecS \vecA^2 \vecU \right) \vece_i \end{smallmatrix} \right).
\end{aligned}
\end{equation}
 where $\vece_1 = \left( \begin{smallmatrix} 1 \\ 0 \end{smallmatrix} \right)$ and $\vece_2 = \left(\begin{smallmatrix} 0 \\ 1 \end{smallmatrix} \right)$.
\\

Now, we want to find some $(\eps, a, \vecX)$ for which the determinant 
\begin{multline}
 \det(\vecX) \: := \:  \biggl| \vecF_1, \vecF_2,   \:  [\vecF_1,\vecF_2], \:  [\vecF_1, [\vecF_1,\vecF_2] ],    \:  [\vecF_2, [\vecF_1,\vecF_2] ], \\
   [\vecF_1, [\vecF_1, [\vecF_1,\vecF_2]]],    \:   [\vecF_2, [\vecF_2, [\vecF_1,\vecF_2] ]]  \biggr|(\vecX) \: \neq 0.
 \end{multline}
As the l.h.s. defines an analytic function of $\vecX$, it will be non-zero almost everywhere. Thus, the Lie algebra generated by $\vecF_1$ and $\vecF_2$ will be maximal (of dimension 7) at almost every $\vecX$, and local controllability will follow from Chow's theorem, see \cite{Jurdjevic97}. 
\\
%
%Finally, the controllability result is based on the following result which is associated with some formal calculation performed by using Maple.
%
%\begin{theorem}
%There exists a function non null $\mathcal{A}$ on $\mathcal{S}$ such that for all $\vecX \in \mathcal{S}$ the asymptotic expansion of $\mbox{det}$ is equal to 
%$$
%\mbox{det}(\vecX)=a^5\,\epsilon^2\, \mathcal{A}(\vecX) + O(a^5\epsilon^2)\,.
%$$
%\end{theorem}
%
%
%\Proof
%
For all $\vecG \in \mbox{Lie}(\vecF_1,\vecF_2)$, let us denote $\vecG^0$, $\vecG^1$ and $\vecG^2$ the zero  order term, the term of order $a$ and the term of order $a \eps$ in the expansion of the vector field $\vecG$ respectively. Thus, 
\begin{equation*}
\vecG = \vecG^0 + \vecG^1 + \vecG^2 + O\left(a^2\right)+O(a\epsilon^2)\,.
\label{eq:asymptotic_expansion_vector_fields}
\end{equation*}
For instance the expansion of the first Lie bracket reads
%$$
%[\vecF_1,\vecF_2] = [\vecF_1^0,\vecF_2^0] +  [\vecF_1^1,\vecF_2^0] +  [\vecF_1^0,\vecF_2^1] + [\vecF_1^2,\vecF_2^0] + [\vecF_1^0,\vecF_2^2] +\left(O\left(a^2\right)+O(a\epsilon^2)\right)\,.
%$$
%It can be readily checked that,
$$
[\vecF_1,\vecF_2] = [\vecF_1,\vecF_2]^0+[\vecF_1,\vecF_2]^1+[\vecF_1,\vecF_2]^2 + O\left(a^2\right)+O(a \epsilon^2)\,.
$$
with 
$$
[\vecF_1,\vecF_2]^0 \: =  \:  [\vecF_1^0,\vecF_2^0] \,,\quad  [\vecF_1,\vecF_2]^1 \: = \:  [\vecF_1^1,\vecF_2^0] +  [\vecF_1^0,\vecF_2^1]\,, \quad  [\vecF_1,\vecF_2]^2 \: = \:  [\vecF_1^2,\vecF_2^0] + [\vecF_1^0,\vecF_2^2].
$$
Note that for all $\vecG \in \mbox{Lie}(\vecF_1,\vecF_2)$, $\vecG^0 + \vecG^1$ is a "flat wall" expansion, first order in $a$. Meanwhile,  $\vecG^2$ is the first  term which takes into account the roughness. 
\\

Without including this extra term, the three-sphere swimmer would not be controllable (see \cite{AlougesGiraldi12}), meaning that the determinant would vanish. We  have notably
\begin{lemma}
For all $\vecG \in \mbox{Lie}(\vecF_1,\vecF_2)\setminus \left\{\vecF_1,\vecF_2\right\}$,  $\vecG^0 = 0$.
\label{lemma_0_order}
\end{lemma}
\Proof
A simple calculation yields
$$
\vecF_1^0(\vecX) = \left(\begin{array}{c} 1\\ 0\\0 \\ 0 \\ \frac{1}{3}\cos(\phi)\sin(\theta)\\\frac{1}{3} \sin(\phi) \sin(\theta)\\ \frac{1}{3} \cos(\theta) \end{array}\right) \,,\quad \quad
\vecF_2^0(\vecX) = \left(\begin{array}{c}0\\ 1\\0 \\ 0 \\ -\frac{1}{3}\cos(\phi)\sin(\theta)\\-\frac{1}{3}\sin(\phi)\sin(\theta)\\ -\frac{1}{3} \cos(\theta) \end{array}\right)\,.
$$
It implies that  $[\vecF_1^0,\vecF_2^0]$ is zero. The lemma is proved. 
\EndProof
As regards the $O(a)$ term, we have 
\begin{lemma} 
Let  $\displaystyle \mbox{Lie}(\vecF_1,\vecF_2)^1  :=  \left\{\vecG^0 + \vecG^1 \,\, s.\,t.\, \vecG \in \mbox{Lie}(\vecF_1,\vecF_2) \right\}$. For all $\vecX$, the dimension of the subspace  $\displaystyle \mbox{Lie}(\vecF_1,\vecF_2)^1(\vecX)$ is less than $5$.
%
%
%For all $\vecG \in \mbox{Lie}(\vecF_1,\vecF_2)$ and for all $\vecX \in \mathcal{S}$, the vector $\vecG^0+\vecG^1(\vecX)$ belongs to the sub-space generated by the vectors $[\vecF_1,\vecF_2]^1(\vecX)$, $[\vecF_2,[\vecF_1,\vecF_2]]^1(\vecX)$, $[\vecF_1,[\vecF_1,\vecF_2]]^1(\vecX)$, i.e., 
%$$
%\vecG^0+\vecG^1(\vecX)\in\mbox{Span}\left([\vecF_1,\vecF_2]^1(\vecX),[\vecF_2,[\vecF_1,\vecF_2]]^1(\vecX),[\vecF_1,[\vecF_1,\vecF_2]]^1(\vecX) \right) 
%$$
%
%The dimension of the Lie algebra generated by $\mbox{Lie}(\vecF_1,\vecF_2)$
%
%
%For all $\vecG \in \mbox{Lie}(\vecF_1,\vecF_2)$, then the leader term $\vecG^0+\vecG^1$ then the fourth component of the first term order in its asymptotic expansion with respect to $a$ and $\epsilon$  small, $\vecG^1\cdot \vece_4$, is null.
\label{lemma_first_order}
\end{lemma}

\Proof
As said above, for all $\vecG \in Lie(\vecF_1,\vecF_2)$,   the sum $\vecG^0 + \vecG^1$ is a $O(a)$ expansion of the "flat wall field", corresponding to the case $h=0$. But in such flat case, symmetries constrain the swimmer within a plane. 
Thus, the associated manifold has at most dimension $5$ ($\xi_1,\xi_2$, two coordinates for the center of the middle ball, one angle). This implies the result. 
\EndProof

\begin{remark}
Since without roughness the swimmer evolves in a plane, it follows that the angle $\phi$ cannot change with time. Consequently, for all $\vecF(\vecX) \in \mbox{Lie}(\vecF_1,\vecF_2)^1(\vecX)$ the fourth component of the vector $\vecF(\vecX)$ is zero.
\label{remark_phi_null}
\end{remark}

\begin{remark}
The lemma  \ref{lemma_first_order} also applies to the vector fields which do not take into account the roughness i.e., the ones which appear in the expansion without $\eps$.
\label{remark_for_all}
\end{remark}

From this, we will get that the non-zero leading term in the expansion of  $\det$  has power $a^5 \epsilon^2$. Theorem \ref{thm:three_sphere}  follows directly from
\begin{proposition} \label{thm:control_final}
In the regime $1 \gg \eps \gg a$,  one can find a surface $h \in C^\infty_c(\R^2)$   and a non-trivial analytic function ${\cal A}$  such that for all $\vecX$ 
$$ \det(\vecX)=a^5\,\epsilon^2\, \mathcal{A}(\vecX) + O(a^6 \eps^2 + a^5\epsilon^3)\,.
$$
\end{proposition}

\Proof
For all  vector $\vecG$, we denote $(\vecG)_{j'}^j \: := \: (\vecG_k)_{j \le k \le j'}$.  Since  $\vecF_i$, $i=1,2$, is of the type $\begin{pmatrix} \vece_i \\ *\\ \vdots \\ * \end{pmatrix}$, we get easily  that 
$$ \det(\vecX) = \left| \vecZ_1, \vecZ_2, \vecZ_3, \vecZ_4, \vecZ_5\right|$$
where 
\begin{equation}
\left\{
\begin{array}{lllll}
\vecZ_1:=([\vecF_1,\vecF_2])_7^3\,,\\
\vecZ_2:=([\vecF_1, [\vecF_1,\vecF_2] ])_7^3\,,\\
\vecZ_3:= ([\vecF_2, [\vecF_1,\vecF_2] ])_7^3\,,\\
\vecZ_4:=([\vecF_1, [\vecF_1, [\vecF_1,\vecF_2]]])_7^3\,,\\
\vecZ_5:= ([\vecF_2, [\vecF_2, [\vecF_1,\vecF_2] ]])_7^3\,.
\end{array}
\right.\,
\label{def_vecZ}
\end{equation}
From  Lemma \ref{lemma_0_order},  $\vecZ_i^0 = 0$ for all $i=1...5$.   Moreover, by Lemma \ref{lemma_first_order}, any determinant of the type 
$$\left| \vecZ^1_{k_1}, \vecZ^1_{k_2}, \vecZ^1_{k_3},  \vecZ^1_{k_4} \right|, \quad k_i \in \{ 1, ..., 5 \} \quad \textrm{is zero.}   
 $$  
Expanding the function $\det$ by 5-linearity, we obtain
$$ \det(\vecX) \: = \: a^5\,\epsilon^2\, \mathcal{A}(\vecX) + O(a^6 \eps^2 + a^5\epsilon^3)\,,
$$
where the function ${\cal A}(\vecX)$ is defined as follows. Let 
$$ {\cal I} \: := \: \left\{ \veck \in \{1,...,5\}^5 \quad \textrm{with} \: k_1 < k_2 \quad \textrm{and} \quad k_3 < k_4 < k_5 \quad \textrm{distinct of $k_1$ and $k_2$}\right\}. $$
We set 
$${\cal A}(\vecX) \: := \: \sum_{\veck \in {\cal I}} \pm \left| \vecZ^2_{k_1}, \:  \vecZ^2_{k_2}, \: \vecZ^1_{k_3}, \:   \vecZ^1_{k_4}, \: \vecZ^1_{k_5} \right|, $$
 where the $\pm$  is the signature of the permutation $i \rightarrow k_i$. 
\\
   
It remains to  prove that there exists $\vecX_0$ such that  $\mathcal{A}(\vecX_0)$ is non-zero. By calling $\vecK^{int}_4$ the function $(\vecs,\vecr,\vecr') \mapsto \frac{\partial}{\partial z}\left(\vecs \mapsto \vecK^0(\vecs,\vecr) \right) \, \frac{\partial}{\partial z}\left(\vecs \mapsto \vecK^0(\vecs,\vecr') \right)$, we have  (see \eqref{K1})
$$
 \vecK_4(\vecr,\vecr') = - \eps \int_{\partial \R_+^3}\,  h(\vecs) \,\vecK^{int}_4(\vecs,\vecr,\vecr')\,  \mbox{d}\vecs\,.
$$

We then define the 3x3 block matrix $\vecA^2_{int}(\vecs)$ through 
$$ (\vecA^2_{int}(\vecs))_{ij} = -6\pi\mu a\, \vecK^{int}_4(\vecs,\vecx_i,\vecx_j), \quad i,j=1...3. 
 $$
By using the linearity of the integral, the vector fields $\vecF_i^2$, $i=1,2$ read 
$$
\vecF^2_i =\left(- \eps \int_{\partial \R_+^3}\,  h(\vecs) \left(\vecF^2_{i,int}(\vecs) \right)\mbox{d}\vecs\right)  \,,
$$
where,
\begin{equation}
\vecF^2_{i,int}(\vecs) =-\left(-(\vecS\vecT)^{-1}\vecS\vecA^2_{int}(\vecs)\vecT(\vecS\vecT)^{-1} \vecS\vecU+ (\vecS\vecT)^{-1} \vecS \vecA^2_{int}(\vecs) \vecU \right)  \vece_i \,.
\end{equation}
Then, denoting 
$$ \vecZ^2_{1,int}(\vecs) \: := \:  [\vecF^2_{1,int}(\vecs),\vecF_2^0] + [\vecF_1^0, \vecF^2_{2,int}(\vecs)]$$
leads to 
\begin{equation}
\vecZ^2_1 = - \eps \int_{\partial \R_+^3}\,  h(\vecs)\, \vecZ^2_{1,int}(\vecs) \,\mbox{d}\vecs\,.
\end{equation}
We can go on with this process and find explicitly functions $\vecZ^2_{i,int}(\vecs)$ for $i=2,\dots,5$ such that  
$$
\forall i\in \{2,\dots,5\}\,, \quad\quad \vecZ_i^2 =- \eps \int_{\partial \R_+^3}\,  h(\vecs) \left(\vecZ^2_{i,int}(\vecs) \right)\mbox{d}\vecs \,.
$$
Finally, 
\begin{multline}
\displaystyle
\mathcal{A}(\vecX) =- \eps^2 \int_{\partial \R_+^3}\,\int_{\partial \R_+^3}\,   h(\vecs)\, h(\vecs')  \\\sum_{\mathcal{I}} \pm \left|\vecZ^2_{k_1,int}(\vecs)\,\,\vecZ^2_{k_2,int}(\vecs')\,\,\vecZ^1_{k_3}\,\,\vecZ^1_{k_4}\,\,\vecZ^1_{k_5}\,(\vecX)\right|\,\mbox{d}\vecs\,\mbox{d}\vecs'.
\end{multline}
We call $\det_{int}$ the function defined by, 
\begin{equation}
\label{eq:detint}
\mbox{det}_{int}:\left(\vecX,\vecs,\vecs'\right) \mapsto \sum_{\mathcal{I}} \pm\left|\vecZ^2_{k_1,int}(\vecs)\,\,\vecZ^2_{k_2,int}(\vecs')\,\,\vecZ^1_{k_3}\,\,\vecZ^1_{k_4}\,\,\vecZ^1_{k_5}\,(\vecX)\right|\,.
\end{equation}

%We call $\det_{int}$ the function $\left(\vecX,\vecs,\vecs'\right) \mapsto \sum_{\mathcal{I}} \pm\left|\vecZ^2_{k_1,int}(\vecs)\,\,\vecZ^2_{k_2,int}(\vecs')\,\,\vecZ^1_{k_3}\,\,\vecZ^1_{k_4}\,\,\vecZ^1_{k_5}\,(\vecX)\right|$.
Clearly, for Theorem \ref{thm:control_final} to hold, it is enough that  there exists $\vecX_0$ and  $\left(\vecs,\vecs'\right) \in (\pa \R^3_+)^2$   such that $\det_{int}(\vecX_0,\vecs,\vecs')$ is not zero for some $\left(\vecs,\vecs'\right)\in \mathbb{R}^4$.  Indeed, we can then adjust the  function $h$ to make the integral non-zero. The calculation of $\det_{int}$ can be carried out using Maple. More precisely, one can derive an equivalent as $z$ goes to infinity, and check that  $\det_{int}(\vecX_0, \cdot, \cdot) \neq 0$ for $\vecX_0 =\left(1,2,\frac{\pi}{3},\frac{\pi}{3},1,2,z\right)$ for $z$ large enough (see appendix \ref{AppendixB} for details). This concludes the proof. 
%\david{z=infty ou z=1 ?}
%
%
%\david{est-ce vraiment utile de mettre une conclusion ? J'ai l'impression que c'est un peu redondant avec l'enonce des resultats}

\section{Conclusions and perspectives}

The aim of this present paper was to examine how the controllability of low Reynolds number artificial swimmers is affected by the presence of a rough wall on a fluid. This study generalizes the one made by F. Alouges and L. Giraldi in \cite{AlougesGiraldi12} which deals with the effect of a plane wall on the controllability of this particular swimmers. 
\\

Firstly, we show Theorem \ref{thm:analyticity}. It deals with the regularity of the dynamics of the swimmers. Indeed, we prove that the equation of motion of such particular swimmers are analytic with respect to the parameters defining the swimmer (radius of the ball, position and length of the arms) and the typical height of roughness of the wall. Then, we deduce Theorem \ref{thm:four_sphere} which claims that the $4$-sphere swimmer remains controllable with the presence of roughness. The proof is based on general arguments which could be used for other models of micro-swimmer. 
% which is expressed as an affine linear control system without drift,
\\

Secondly, Theorem \ref{thm:three_sphere} examines the controllability of the Three-sphere swimmer in the presence of a rough wall. More precisely, we show that there exists a roughness such that the swimmer can locally reach any direction. We recall that the previous studies made on the 3-sphere swimmer allow to show that it can reach only one direction (see \cite{AlougesDeSimone10} when it evolves in a whole space and three directions with the presence of a plane wall (see \cite{AlougesGiraldi12}). In our case, the roughness leads to break the symmetry of the system "fluids-swimmer". As a result, it allows the swimmer to reach any direction. The proof is an in-depth study which associates several tools both in hydrodynamics and control theory. The general "idea" emphasizes here is the fact that in the real life all the micro-organism, regardless how symmetric it is, can move in any direction.  
%based on several steps which associate the regularity of the dynamics which leads to apply Chow Theorem and also formal calculation
\\

The quantitative approach to this question together with the complete understanding in a view of controllability of underlying systems is far beyond reach and thus still under progress as in a another direction, the consideration of an confined environment, e.g. when the fluid is bounded.   
Future work will also explore which are the directions easier to reach than the others by varying the rough wall.  

\newpage
\bibliographystyle{plain}
\bibliography{ReferencesRugosite}
\newpage
\appendix
\section{A well-posedness result for the Stokes system}
We show here the well-posedness of the inhomogeneous Stokes system involved in the proof of Theorem \ref{thm:analyticity}. We refer to this proof for notations, and shall drop here all bars for brevity. What we want to show is  
\begin{theorem}
Let $(F,G,V_1,...,V_N)$ given in $\ds {\cal V}' \times L^2({\cal F}) \times \prod_l H^{1/2}(B_l)$. There exists a unique solution $(V,Q)$ in ${\cal V}_0 \times L^2({\cal F})$ of 
\begin{align*}
-\Delta V  + \na Q & = F \quad  \: \mbox{in} \: {\cal F}, \\
 \Div V & =  G \quad  \: \mbox{in} \: {\cal F}, \\ 
 V & = V_l   \quad \mbox{at} \: \pa B_l, \quad l=1...N,
\end{align*}
\end{theorem}
We recall that the space ${\cal V}_0$ encodes the additional homogeneous Dirichlet condition at $\pa {\cal O}$. 

\medskip
\noindent
{\em Proof of the Theorem}. The theorem follows from 
\begin{proposition} \label{propdivergence}
For all $(G,V_1,...,V_n)$   given in $\ds L^2({\cal F}) \times \prod_l H^{1/2}(B_l)$, there exists a field $W \in {\cal V}_0$ satisfying
\begin{equation} \label{divdir}
\Div W =  G \quad  \: \mbox{in} \: {\cal F}, \quad 
\quad W  = V_l   \quad \mbox{at} \: \pa B_l, \quad l=1...N. 
\end{equation}
together with the estimate: $\| \na W \|_{L^2} \le C \left( \|Â G \|_{L^2} + \sum_{l=1}^N  \| V_l \|_{H^{1/2}} \right)$. 
\end{proposition}
This proposition will be proved below. Let us explain how it implies the theorem. First, considering  $V' := V - W$, and $F' := F + \Delta W$  one can come down to the homogeneous case  $G=0$ and $V_l = 0$ for all $l$.  
The homogeneous case can then be solved by a standard application of Lax-Milgram theorem. More precisely, defining 
$$   {\cal V}_{hom,div} \: := \: \left\{V' \in {\cal V}_0, \quad V'\vert_{\pa B_l} = 0 \:\: \forall l, \quad \div V' = 0 \right\} $$    
one can show easily that there is a unique $V' \in  {\cal V}_{hom,div}$ satisfying 
$$ \int \na V' \cdot \na \phi \: = \:  < F, \phi >, \quad \forall \phi \in  {\cal V}_{hom,div}. $$
We recall that the condition $V/(1+|\vecx|) \in L^2({\cal F})$ in the definition of ${\cal V}_0$ is related to the Hardy inequality.  

\medskip
By standard arguments, one then recovers a pressure field $Q \in L^2_{loc}({\cal F})$ so that the Stokes equation $-\Delta V' + \na Q = F'$ holds. Eventually, to show that we can take $Q$ in $L^2({\cal F})$, we invoke \cite[Theorem 3.5.3, page 217]{Galdi11}: it is enough that for all $G \in L^2({\cal F})$, the problem 
$$  \div W = G \:\:  \mbox{in} \: {\cal F}, \quad W\vert_{\pa {\cal F}} = 0  $$
has a solution $W \in {\cal V}_0$ with  $\| \na W \|_{L^2} \le C \| G \|_{L^2}$. This is a special case of Proposition \ref{propdivergence}. This ends the proof. 

\bigskip
\noindent
{\em Proof of the Proposition}. Again, we single out the key ingredient in a 
\begin{lemma}
Given $G \in L^2({\cal O})$, there exists a field  $W \in {\cal V}_0$ 
such that 
$$ \Div W  = G, \quad  W\vert_{\pa {\cal O}} = 0 \quad   \| \na W \|_{L^2} \: \le \:  C \, \| G \|_{L^2}.  $$  
\end{lemma}
Note that this lemma is only about the domain ${\cal O}$, that is without the balls. Let us postpone its proof, and show how it implies the existence of a $W$ satisfying \eqref{divdir}. 
\begin{itemize}
\item First step:  we lift the boundary data $V_l$. One can find $\overline{W} \in H^1({\cal F})$ compactly supported near the balls, such that $\overline{W} = V_l$ at $\pa B_l$. Up to replace $W$ by $W - \overline{W}$ and $G$ by   $G - \Div \overline{W}$, we can assume $V_l = 0$ for all $l$. 
\item Second step (assuming now $V_l = 0$ for all $l$): we extend $G$ by $0$ in the balls and apply the Lemma: it provides a $\tilde W$ satisfying  $ \Div  \tilde W  =  G, \quad   \tilde W\vert_{\pa {\cal O}} = 0$.  However, the boundary data at the balls is non-zero: $\tilde W\vert_{\pa B_l} \neq 0$. 
\item   Third step: we correct this non-zero boundary data. We observe that 
$$\int_{\pa B_l} \tilde W \cdot n \, ds= 0 = \int_{B_l} \Div \tilde W = 0,$$ 
as $G$ was extended by  zero inside the balls. Thanks to this "compatibility" condition, we can use a standard result of Bogovskii, see \cite[Exercice III.3.5, p176]{}: for all $l$, there exists a field $W_l$ defined over the annulus $\{a < |x - x_l| < a+\eta\}$, satisfying 
$$ \div W_l = 0, \quad W_l\vert_{\pa B_l} = - \tilde W\vert_{\pa B_l}, \quad W_l\vert_{\{|x-x_l| = a + \eta\}} = 0.   $$          
We take $\eta$ small enough so that the annuli do not intersect. Then,  we extend the $W_l$'s by $0$ outside the annuli and set $W := \tilde W + \sum W_l$. This new field $W$ satisfies \eqref{divdir}, as expected.  
\end{itemize}

\noindent
{\em Proof of the Lemma}. In the case where $h=cst$, that is for a flat half-space, the result is  classical: {\it cf} \cite[Corollary 4.3.1, p261]{}.  In particular, if the support of $G$ is included in $\{x_3 >  \sup |h|\}$, the problem is solved: one can take the solution $W$ of    
$$ \Div W  = G \quad \mbox{for} \quad  x_3 > \sup |h|, \quad  W\vert_{\{x_3 = \sup |h|\}} = 0   $$  
and extend it by zero below $\{x_3 = \sup |h|\}$. 
\\

For general $G$, we can decompose $G \: = \:  G \, 1_{ \{ x_3 > \sup |h|}\} \: + \:  G \, 1_{ \{ x_3 < \sup |h|\}}$, and handle the first part as previously. In other words, {\em it remains to consider the case where $G$ is compactly supported in $x_3$.}  From there, we proceed in three steps: 
\begin{itemize}
\item Step 1. Let $R$ such that $G = 0$ for $x_3 \ge R$.  We introduce $W^1 := \na \psi \, 1_{\{x_3 < R\}}$ where $\psi$ satisfies  
$$ \Delta \psi = G   \quad \mbox{for} \: \eps h < x_3 < R, \quad \pa_n \psi\vert_{\pa {\cal O}} = 0, \quad \psi\vert_{x_3 =R} = 0. $$
This Poisson equation has a unique solution in $H^2(\{ \eps h < x_3 < R \})$:  note that Poincar\'e inequality applies thanks to the Dirichlet condition at $x_3 = R$. Hence, $W^1$ satisfies $\Div W^1 = G$ in the strip $\{ \eps h < x_3 < R \}$, and also trivially in the half-space $\{x_3  > R \}$. However, two problems remain: the normal component of $W^1$ jumps at $x_3 = R$, and it has non-zero boundary data at $\{ x_3 = \eps h \}$. 
\item Step 2. Correction of the jump at $x_3 = R$. We just introduce the field $W^2 := \tilde W \,  1_{\{x_3 > R\}}$, where $\tilde W$ satisfies 
\begin{multline*} \Div \tilde W = 0  \quad \mbox{for} \quad  x_3 > R,   \quad \tilde W\vert_{\{x_3 = R\}} = \na \Psi\vert_{x_3 = R}, \\
 \| \na \tilde W \|_{L^2} \: \le \:  C \| \na \psi \|_{H^{1/2}(\{x_3 = R \})} \:  (\le \: C \,   \|G \|_{L^2}).  
 \end{multline*}
The existence of such $\tilde W$ is classical, see \cite[Theorem IV.3.3]{Galdi11}. 
\item Step 3. Correction of the boundary data. Thanks to the Neumann condition on $\psi$, we have $W^1 \cdot n\vert_{\pa {\cal O}} = 0$. We introduce some partition of unity $(\chi_k = \chi_k(x_1,x_2))_{k \in \Z^2}$ associated to a covering of $\R^2$  by rectangles $R_k$. More precisely, we assume that  the lengths of $R_k$ are   uniformly bounded in $k$, and that the $C^1$ norms of $\chi^k$ are uniformly bounded in $k$ (we leave the construction of examples to the reader).     
Thanks to the tangency condition on   $W^1$, we can apply the Bogovskii's result seen above on slices  $S_k \: := \:  \{ (x_1,x_2) \in R_k,  \quad  \: \eps h(x_1,x_2) < x_3 < R \}, \: k \in \Z^2$. Hence, there exists some $W_k \in H^1(S_k)$ such that 
\begin{multline*}
 \Div W_k = 0 \quad \mbox{in} \quad S_k, \quad W_k = -\chi_k W^1 \quad  \mbox{at} \quad  \pa S_k \cap \pa {\cal O},   \\
  W_k = 0   \quad \mbox{at} \quad   \pa S_k \setminus {\cal O},
  \end{multline*} 
 and $\| \na W_k \|_{L^2} \: \le \: C \, \|\chi_k W^1 \|_{H^{1/2}(\pa {\cal O})}$. Extending all $W_k$'s by $0$ outside $S_k$, and setting $W^3 \: := \: \sum_{k \in \Z^2}   W_k$, we find that 
 \begin{multline*} 
 \div W^3 = 0 \quad \mbox{in} \quad  { \cal O }, \quad W^3\vert_{\pa {\cal O}} \: = - W^1\vert_{\pa {\cal O}}, \\ 
 ||\na W^3\|_{L^2({\cal O})} \: \le \: C \| W^1 \|_{H^{1/2}(\pa {\cal O})} \: (\le C \|W^1 \|_{H^1({\cal O})}).
 \end{multline*}
\end{itemize}
Finally, $W = W^1 + W^2 + W^3$ fulfills all requirements, which concludes the proof of the lemma.

\section{Formal expressions}
 \label{AppendixB}
We express here the requisite formal expression of the vector fields for the calculus of the $\det_{int}$ at the point $\vecX_0 =\left(1,2,\frac{\pi}{3},\frac{\pi}{3},1,2,z\right)$.
\\

First at all, we have used the software MAPLE to symbolically compute The vector fields $\vecF_1$ and $\vecF_2$  by using the formula \eqref{eq:vectorfields}. Then, we deduce the expression of every  vector which belongs to the set $\mathcal{F}_{cal}:=\left\{\vecZ_k^j \quad \textrm{s.t.}\quad k=1,\cdots,5 \,\, j=1,2 \right\}$ defined in \eqref{def_vecZ}. In the following, we express the first asymptotic terms when $z$ goes to infinity of the vector fields which belong to $\mathcal{F}_{cal}$ at $\vecX_0$. The asymptotic expression of the determinant $\det_{int}$, defined in \eqref{eq:detint}, is deduced.

\begin{itemize}
\item The expansion of vector fields $\vecZ^1_1$ and $\vecZ_1^2$ are expressed by, 
$$
\vecZ_1^1(\vecX_0) = \left(\begin{smallmatrix}
\frac{21627}{57344}\frac{\sqrt{3}}{z^4} - \frac{237897}{802816}\frac{\sqrt{3}}{z^5} - \frac{56965095}{25690112}\frac{\sqrt{3}}{z^6} - \frac{29418201}{25690112}\frac{\sqrt{3}}{z^7}+\frac{141}{25088}\frac{\sqrt{3}}{z^8}\\
 0\\ 
\frac{41}{432}\sqrt{3} - \frac{7209}{229376}\frac{\sqrt(3)}{z^4} -\frac{923829}{3211264}\frac{\sqrt(3)}{z^5}+\frac{45738445}{102760448}\frac{\sqrt{3}}{z^6} -\frac{83758845}{102760448}\frac{\sqrt(3)}{z^7}-\frac{47}{100352}\frac{\sqrt{3}}{z^8}\\
\frac{41}{144} -\frac{21627}{229376 z^4} - \frac{2771487}{3211264 z^5} + \frac{137215335}{102760448 z^6} - \frac{251276535}{102760448 z^7} - \frac{141}{100352 z^8}\\
\frac{41}{216}+\frac{21627}{114688 z^4} - \frac{927033}{1605632 z^5} - \frac{38588135}{51380224 z^6} + \frac{191091795}{51380224 z^7} - \frac{61}{100352 z^8}
\end{smallmatrix}\right)\,,
$$

%$$
%\vecZ_1^2(\vecX_0,s,s') = \left(\begin{smallmatrix} 
%\frac{2187}{448}\frac{\left(\frac{1}{2}+\frac{1}{2}\sqrt{3}\right)\left(1-\frac{1}{4}\sqrt{3}\right)}{\pi z^7} 
%-\frac{27}{784}\frac{\sqrt{3}\left(\frac{19593}{64}v + \frac{6531}{64} s -\frac{567}{64}\sqrt{3}s' - \frac{189}{64}\sqrt{3}s - \frac{927}{128}\sqrt{3} - \frac{63609}{64}\right)}{\pi z^8}\\
%-\frac{2187}{448}\frac{-\frac{1}{2}\sqrt(3)+\frac{1}{2}}{\pi z^7} 
%+ \frac{27}{784}\frac{\sqrt{3}\left(-\frac{6531}{16}s'+\frac{7665}{16}s - \frac{189}{16}\sqrt{3}s' - \frac{189}{16}\sqrt{3}s - \frac{2997}{32}\sqrt{3} + \frac{5727}{8}\right)}{\pi z^8}\\
%\frac{81}{448} \frac{\left(-\frac{39}{32}\sqrt(3)+\frac{81}{32}\right)\sqrt(3)}{\pi z^7} 
%- \frac{9}{196} \frac{\frac{8505}{1024}s' + \frac{1587915}{1024}s+\frac{58779}{1024}\sqrt{3}s' - \frac{3087}{4}s^2 + \frac{230769}{1024}\sqrt{3}s -\frac{406647}{1024}\sqrt{3} - \frac{2126079}{2048}}{\pi z^8}\\
%-\frac{243}{448}\frac{-\frac{9}{32}\sqrt{3} + \frac{63}{32}}{\pi z^7} + \frac{9}{196}\frac{-\frac{4011777}{1024}s' - \frac{72387}{1024}s + \frac{567}{1024}\sqrt{3}s' + \frac{3087}{4}s'^2 + \frac{1701}{1024}\sqrt(3)s+\frac{54351}{2048}\sqrt{3} + \frac{5792805}{1024}}{\pi z^8}\\
%\frac{729}{896}\frac{\left(\frac{1}{2} + \frac{1}{2}\sqrt{3}\right)\left(1 - \frac{1}{4}\sqrt{3}\right)\sqrt{3}}{\pi z^7} - \frac{27}{1568}\frac{\frac{19593}{64}s' + \frac{6531}{64}s - \frac{567}{64}\sqrt{3}s' - \frac{189}{64}\sqrt{3}s - \frac{927}{128}\sqrt{3} - \frac{63609}{64}}{\pi z^8}
%\end{smallmatrix}\right)
%$$

$$
\begin{array}{ll}
\vecZ_1^2(\vecX_0,s,s') =& \left(\begin{smallmatrix} 
\frac{2187}{448}\frac{\left(\frac{1}{2}+\frac{1}{2}\sqrt{3}\right)\left(1-\frac{1}{4}\sqrt{3}\right)}{\pi z^7} \\
-\frac{2187}{448}\frac{-\frac{1}{2}\sqrt(3)+\frac{1}{2}}{\pi z^7} \\
\frac{81}{448} \frac{\left(-\frac{39}{32}\sqrt(3)+\frac{81}{32}\right)\sqrt(3)}{\pi z^7} \\
-\frac{243}{448}\frac{-\frac{9}{32}\sqrt{3} + \frac{63}{32}}{\pi z^7}\\
\frac{729}{896}\frac{\left(\frac{1}{2} + \frac{1}{2}\sqrt{3}\right)\left(1 - \frac{1}{4}\sqrt{3}\right)\sqrt{3}}{\pi z^7}
\end{smallmatrix}\right)\quad+\\
& \left(\begin{smallmatrix} 
-\frac{27}{784}\frac{\sqrt{3}\left(\frac{19593}{64}v + \frac{6531}{64} s -\frac{567}{64}\sqrt{3}s' - \frac{189}{64}\sqrt{3}s - \frac{927}{128}\sqrt{3} - \frac{63609}{64}\right)}{\pi z^8}\\
\frac{27}{784}\frac{\sqrt{3}\left(-\frac{6531}{16}s'+\frac{7665}{16}s - \frac{189}{16}\sqrt{3}s' - \frac{189}{16}\sqrt{3}s - \frac{2997}{32}\sqrt{3} + \frac{5727}{8}\right)}{\pi z^8}\\
- \frac{9}{196} \frac{\frac{8505}{1024}s' + \frac{1587915}{1024}s+\frac{58779}{1024}\sqrt{3}s' - \frac{3087}{4}s^2 + \frac{230769}{1024}\sqrt{3}s -\frac{406647}{1024}\sqrt{3} - \frac{2126079}{2048}}{\pi z^8}\\
 \frac{9}{196}\frac{-\frac{4011777}{1024}s' - \frac{72387}{1024}s + \frac{567}{1024}\sqrt{3}s' + \frac{3087}{4}s'^2 + \frac{1701}{1024}\sqrt(3)s+\frac{54351}{2048}\sqrt{3} + \frac{5792805}{1024}}{\pi z^8}\\
 - \frac{27}{1568}\frac{\frac{19593}{64}s' + \frac{6531}{64}s - \frac{567}{64}\sqrt{3}s' - \frac{189}{64}\sqrt{3}s - \frac{927}{128}\sqrt{3} - \frac{63609}{64}}{\pi z^8}
\end{smallmatrix}\right)\\
& +\quad O(\frac{1}{z^9})\,.
\end{array}
$$
\item The expansion of vector fields $\vecZ^1_2$ and $\vecZ_2^2$ are expressed by,

$$
\vecZ_2^1(\vecX_0) = \left(\begin{smallmatrix}
\frac{7209}{25088}\frac{\sqrt{3}}{z^4} - \frac{2097819}{11239424}\frac{\sqrt{3}}{z^5} - \frac{138036945}{44957696}\frac{\sqrt{3}}{z^6} + \frac{166111299}{89915392}\frac{\sqrt{3}}{z^7} + \frac{477316227}{359661568}\frac{\sqrt{3}}{z^8} \\
0\\
-\frac{13}{81}\sqrt(3) - \frac{2403}{100352}\frac{\sqrt{3}}{z^4} -\frac{18025783}{44957696}\frac{\sqrt{3}}{z^5} + \frac{325287505}{539492352}\frac{\sqrt{3}}{z^6} - \frac{852035953}{359661568}\frac{\sqrt(3)}{z^7} + \frac{1369121487}{1438646272}\frac{\sqrt{3}}{z^8}\\
-\frac{13}{27} - \frac{7209}{100352 z^4} - \frac{54077349}{44957696 z^5} + \frac{325287505}{179830784 z^6} - \frac{2556107859}{359661568 z^7} + \frac{4107364461}{1438646272 z^8}\\
-\frac{26}{81} + \frac{7209}{50176 z^4} - \frac{14961691}{22478848 z^5} - \frac{285472115}{269746176 z^6} + \frac{2178996509}{179830784 z^7} - \frac{3139098785}{719323136 z^8}
\end{smallmatrix}\right)\,,
$$

%$$
%\vecZ_2^2(\vecX_0,s,s') = \left(\begin{smallmatrix} 
%\frac{729}{196} \frac{\left(\frac{1}{2} + \frac{1}{2}\sqrt{3}\right)\left(1-\frac{1}{4}\sqrt{3}\right)}{\pi z^7} 
%- \frac{9}{5488} \frac{\sqrt{3} \left( \frac{11403}{8}s' + \frac{3801}{8}s - \frac{567}{4}\sqrt{3}s' - \frac{189}{4}\sqrt{3}s + \frac{99591}{128}\sqrt{3} + \frac{26481}{64}\right)}{\pi z^8}\\
%-\frac{729}{196} \frac{-\frac{1}{2}\sqrt{3} + \frac{1}{2}}{\pi z^7} + \frac{9}{5488}\sqrt{3} \frac{\left( -4161 - \frac{3801}{2}s' + \frac{6069}{2}s -189\sqrt{3}s' -189\sqrt{3}s+\frac{69765}{32}\sqrt{3} \right)}{\pi z^8}\\
%\frac{9}{3136}\sqrt{3} \frac{\left( \frac{243}{2} - \frac{117}{2}\sqrt{3}\right)}{\pi z^7} -\frac{3}{1372} \frac{\frac{8505}{64}s' + \frac{698859}{64} s + \frac{34209}{128}\sqrt{3}s' - \frac{21609}{4}s^2 - \frac{307125}{128}\sqrt{3}s+\frac{2517615}{1024}\sqrt{3} - \frac{23258601}{2048}}{\pi z^8}\\
%-\frac{27}{3136} \frac{-\frac{27}{2}\sqrt{3} + \frac{189}{2}}{\pi z^7} + \frac{3}{1372} \frac{\left(-\frac{1994643}{128}s' - \frac{61425}{128} s + \frac{567}{64}\sqrt{3}s' + \frac{21609}{4}s'^2 + \frac{1701}{64}\sqrt{3} s + \frac{54297}{2048}\sqrt{3} + \frac{26018883}{1024}\right)}{\pi z^8}\\
%\frac{243}{392} \frac{\left(\frac{1}{2} + \frac{1}{2}\sqrt{3}\right) \left(1 - \frac{1}{4}\sqrt{3}\right)\sqrt{3}}{\pi z^7} 
%+ \frac{9}{10976} \frac{\left(-\frac{11403}{8}s' - \frac{3801}{8}s  + \frac{567}{4}\sqrt{3}s' + \frac{189}{4}\sqrt{3} s - \frac{26481}{64} - \frac{99591}{128}\sqrt{3}\right)}{\pi z^8}\\
%\end{smallmatrix}\right)
%$$

$$
\begin{array}{lll}
\vecZ_2^2(\vecX_0,s,s') =& \left(\begin{smallmatrix} 
\frac{729}{196} \frac{\left(\frac{1}{2} + \frac{1}{2}\sqrt{3}\right)\left(1-\frac{1}{4}\sqrt{3}\right)}{\pi z^7} \\
-\frac{729}{196} \frac{-\frac{1}{2}\sqrt{3} + \frac{1}{2}}{\pi z^7} \\
\frac{9}{3136}\sqrt{3} \frac{\left( \frac{243}{2} - \frac{117}{2}\sqrt{3}\right)}{\pi z^7} \\
-\frac{27}{3136} \frac{-\frac{27}{2}\sqrt{3} + \frac{189}{2}}{\pi z^7} \\
\frac{243}{392} \frac{\left(\frac{1}{2} + \frac{1}{2}\sqrt{3}\right) \left(1 - \frac{1}{4}\sqrt{3}\right)\sqrt{3}}{\pi z^7} \\
\end{smallmatrix}\right) \quad + \\
& \left(\begin{smallmatrix} 
- \frac{9}{5488} \frac{\sqrt{3} \left( \frac{11403}{8}s' + \frac{3801}{8}s - \frac{567}{4}\sqrt{3}s' - \frac{189}{4}\sqrt{3}s + \frac{99591}{128}\sqrt{3} + \frac{26481}{64}\right)}{\pi z^8}\\
+ \frac{9}{5488}\sqrt{3} \frac{\left( -4161 - \frac{3801}{2}s' + \frac{6069}{2}s -189\sqrt{3}s' -189\sqrt{3}s+\frac{69765}{32}\sqrt{3} \right)}{\pi z^8}\\
 -\frac{3}{1372} \frac{\frac{8505}{64}s' + \frac{698859}{64} s + \frac{34209}{128}\sqrt{3}s' - \frac{21609}{4}s^2 - \frac{307125}{128}\sqrt{3}s+\frac{2517615}{1024}\sqrt{3} - \frac{23258601}{2048}}{\pi z^8}\\
+ \frac{3}{1372} \frac{\left(-\frac{1994643}{128}s' - \frac{61425}{128} s + \frac{567}{64}\sqrt{3}s' + \frac{21609}{4}s'^2 + \frac{1701}{64}\sqrt{3} s + \frac{54297}{2048}\sqrt{3} + \frac{26018883}{1024}\right)}{\pi z^8}\\
+ \frac{9}{10976} \frac{\left(-\frac{11403}{8}s' - \frac{3801}{8}s  + \frac{567}{4}\sqrt{3}s' + \frac{189}{4}\sqrt{3} s - \frac{26481}{64} - \frac{99591}{128}\sqrt{3}\right)}{\pi z^8}\\
\end{smallmatrix}\right)\\
& +\quad O(\frac{1}{z^9})\,.
\end{array}
$$

\item The expansion of vector fields $\vecZ^1_3$ and $\vecZ_3^2$ are expressed by, 

$$
\vecZ_3^1(\vecX_0) = \left(\begin{smallmatrix}
\frac{36045}{802816}\frac{\sqrt{3}}{z^4} - \frac{1737369}{22478848}\frac{sqrt{3}}{z^5}- \frac{688526655}{359661568}\frac{\sqrt{3}}{z^6} - \frac{1554687891}{359661568}\frac{\sqrt{3}}{z^7} - \frac{473853315}{359661568}\frac{\sqrt{3}}{z^8}\\
0\\
-\frac{19}{1296}\sqrt{3} - \frac{12015}{3211264}\frac{\sqrt{3}}{z^4} - \frac{22827197}{89915392}\frac{\sqrt{3}}{z^5} + \frac{1624779455}{4315938816}\frac{\sqrt{3}}{z^6} - \frac{975961023}{1438646272}\frac{\sqrt(3)}{z^7} -\frac{1370275791}{1438646272}\frac{\sqrt{3}}{z^8}\\
-\frac{19}{432} - \frac{36045}{3211264 z^4} - \frac{68481591}{89915392 z^5} + \frac{1624779455}{1438646272 z^6} - \frac{2927883069}{1438646272 z^7} - \frac{4110827373}{1438646272 z^8}\\
-\frac{19}{648}+\frac{36045}{1605632 z^4} -\frac{17817209}{44957696 z^5} -\frac{1422386365}{2157969408 z^6} + \frac{2266030269}{719323136 z^7} + \frac{3127809953}{719323136 z^8}\\
\end{smallmatrix}\right)\,,
$$

$$
\begin{array}{lll}
\vecZ_3^2(\vecX_0,s,s') =& \left(\begin{smallmatrix}
\frac{3645}{6272} \frac{\left(\frac{1}{2} + \frac{1}{2}\sqrt{3}\right)\left(1- \frac{1}{4}\sqrt{3}\right)}{\pi z^7} \\
-\frac{3645}{6272}\frac{\left(-\frac{1}{2}\sqrt{3} + \frac{1}{2}\right)}{\pi z^7} \\
\frac{9}{3136}\sqrt{3}\frac{\left(-\frac{585}{64}\sqrt{3}+\frac{1215}{64}\right)}{\pi z^7}\\ 
-\frac{27}{3136}\frac{\left(\frac{945}{64} - \frac{135}{64}\sqrt{3}\right)}{\pi z^7} \\
\frac{1215}{12544}\frac{\left(\frac{1}{2}+\frac{1}{2}\sqrt{3}\right)\left(1-\frac{1}{4}\sqrt{3}\right)\sqrt{3}}{\pi z^7}\\
\end{smallmatrix}\right)  \quad +\\
&\left(\begin{smallmatrix}
- \frac{9}{5488}\sqrt{3}\frac{\left(\frac{320229}{128}s'+\frac{106743}{128}s-\frac{2835}{128}\sqrt{3}s' - \frac{945}{128}\sqrt{3}s - \frac{48927}{64}\sqrt{3} - \frac{214935}{16}\right)}{\pi z^8}\\
 - \frac{9}{5488}\sqrt{3}\frac{\left(\frac{106743}{32}s' - \frac{112413}{32} s + \frac{945}{32}\sqrt{3}s' - \frac{213405}{16} + \frac{945}{32}\sqrt{3}s + \frac{53283}{16}\sqrt{3}\right)}{\pi z^8}\\
 - \frac{3}{1372} 
\frac{\left(\frac{42525}{2048}s' + \frac{22164471}{2048}s + \frac{960687}{2048}\sqrt{3}s' - \frac{21609}{4}s^2+\frac{11452077}{2048}\sqrt{3}s - \frac{2054187}{256}\sqrt{3} -\frac{10591731}{1024}\right)}{\pi z^8}\\ 
 - \frac{3}{1372}\frac{\left(\frac{80736957}{2048}s'+\frac{1028727}{2048}s - \frac{2835}{2048}\sqrt{3}s' - \frac{21609}{4}s'^2- \frac{8505}{2048}\sqrt{3}s - \frac{19271313}{256} - \frac{677835}{1024}\sqrt{3}\right)}{\pi z^8}\\
-\frac{9}{10976} \frac{\left(\frac{320229}{128}s'+\frac{106743}{128}s - \frac{2835}{128}\sqrt{3}s'-\frac{945}{128}\sqrt{3}s-\frac{48927}{64}\sqrt{3}-\frac{214935}{16}\right)}{\pi z^8}\\
\end{smallmatrix}\right)\\
& +\quad O(\frac{1}{z^9})\,.
\end{array}
$$
\item The expansion of vector fields $\vecZ^1_4$ and $\vecZ_4^2$ are expressed by,

$$
\vecZ_4^1(\vecX_0) = \left(\begin{smallmatrix}
-\frac{64881}{351232}\frac{\sqrt{3}}{z^4}+\frac{69687}{307328}\frac{\sqrt{3}}{z^5}-\frac{1076724525}{629407744}\frac{\sqrt{3}}{z^6}+\frac{5383125903}{629407744}\frac{\sqrt{3}}{z^7}-\frac{4471078809}{2517630976}\frac{\sqrt{3}}{z^8}\\
0\\
\frac{40}{81}\sqrt{3}+\frac{21627}{1404928}\frac{\sqrt{3}}{z^4} - \frac{2233885}{9834496}\frac{\sqrt{3}}{z^5} +\frac{2387478445}{7552892928}\frac{\sqrt{3}}{z^6} - \frac{28725839983}{7552892928}\frac{\sqrt{3}}{z^7}+\frac{50480747763}{10070523904}\frac{\sqrt{3}}{z^8}\\
\frac{40}{27}+\frac{64881}{1404928 z^4}-\frac{6701655}{9834496 z^5}+\frac{2387478445}{2517630976 z^6} - \frac{28725839983}{2517630976 z^7}+\frac{151442243289}{10070523904 z^8}\\
\frac{80}{81}-\frac{64881}{702464 z^4}-\frac{424745}{2458624 z^5}-\frac{2329702535}{3776446464 z^6}+\frac{71414731349}{3776446464 z^7}-\frac{138596989329}{5035261952 z^8}
\end{smallmatrix}\right)\,,
$$

%$$
%\vecZ_4^2(\vecX_0,s,s') = \left(\begin{smallmatrix}
% -\frac{6561}{2744}\frac{\left(\frac{1}{2}+\frac{1}{2}\sqrt{3}\right) \left(1-\frac{1}{4}\sqrt{3}\right)}{\pi z^7} - \frac{27}{38416} \left(\frac{1}{2}+\frac{1}{2}\sqrt{3}\right) \left(1-\frac{1}{4}\sqrt{3}\right) \frac{\left(5103 s'+1701 s+74197/4-\frac{1839}{4}\sqrt{3}\right)}{\pi z^8}\\
%  \frac{6561}{2744}\frac{\left(-\frac{1}{2}\sqrt{3}+\frac{1}{2}\right)}{\pi z^7} + \frac{27}{38416}\frac{\left(\frac{1701}{2}s'+\frac{1701}{2}s - \frac{5103}{2}\sqrt{3} s'+\frac{1701}{2}\sqrt{3}s-\frac{22411}{2}\sqrt{3}+\frac{53465}{4}\right)}{\pi z^8}\\
% \frac{6561}{2744} \frac{\left(-\frac{1}{2}\sqrt{3}+\frac{1}{2}\right)}{\pi z^7} + \frac{27}{38416} \frac{\left(\frac{1701}{2}s'+\frac{1701}{2}s-\frac{5103}{2}\sqrt{3}s'+\frac{1701}{2}\sqrt{3}s-\frac{22411}{2}\sqrt{3}+\frac{53465}{4}\right)}{\pi z^8}\\
% -\frac{27}{21952}\frac{\left(\frac{729}{4}\sqrt{3}-\frac{1053}{4}\right)}{\pi z^7}+\frac{1}{19208}\sqrt{3} \frac{\left(-\frac{413343}{64}s'+\frac{7490259}{64}s+\frac{76545}{64}\sqrt{3}s'+\frac{66339}{64}\sqrt{3}s+\frac{7226181}{128}\sqrt{3}-\frac{4527333}{32}\right)}{\pi z^8}\\
% -\frac{27}{21952}\sqrt{3} \frac{\left(-\frac{567}{4}\sqrt{3}+\frac{243}{4}\right)}{\pi z^7} - \frac{3}{19208} \frac{\left(-\frac{7582113}{64}s' + \frac{15309}{64}s+/frac{5103}{64}\sqrt{3}s'+\frac{15309}{64}\sqrt{3}s+\frac{295803}{128} \sqrt{3}+\frac{2697651}{32}\right)}{\pi z^8}
%\end{smallmatrix}\right)
%$$

$$
\begin{array}{lll}
\vecZ_4^2(\vecX_0,s,s') = &\left(\begin{smallmatrix}
 -\frac{6561}{2744}\frac{\left(\frac{1}{2}+\frac{1}{2}\sqrt{3}\right) \left(1-\frac{1}{4}\sqrt{3}\right)}{\pi z^7} \\
  \frac{6561}{2744}\frac{\left(-\frac{1}{2}\sqrt{3}+\frac{1}{2}\right)}{\pi z^7} \\
 \frac{6561}{2744} \frac{\left(-\frac{1}{2}\sqrt{3}+\frac{1}{2}\right)}{\pi z^7} \\
 -\frac{27}{21952}\frac{\left(\frac{729}{4}\sqrt{3}-\frac{1053}{4}\right)}{\pi z^7}\\
 -\frac{27}{21952}\sqrt{3} \frac{\left(-\frac{567}{4}\sqrt{3}+\frac{243}{4}\right)}{\pi z^7}
 \end{smallmatrix}\right)\quad + \\
& \left(\begin{smallmatrix}
 - \frac{27}{38416} \left(\frac{1}{2}+\frac{1}{2}\sqrt{3}\right) \left(1-\frac{1}{4}\sqrt{3}\right) \frac{\left(5103 s'+1701 s+74197/4-\frac{1839}{4}\sqrt{3}\right)}{\pi z^8}\\
  + \frac{27}{38416}\frac{\left(\frac{1701}{2}s'+\frac{1701}{2}s - \frac{5103}{2}\sqrt{3} s'+\frac{1701}{2}\sqrt{3}s-\frac{22411}{2}\sqrt{3}+\frac{53465}{4}\right)}{\pi z^8}\\
 + \frac{27}{38416} \frac{\left(\frac{1701}{2}s'+\frac{1701}{2}s-\frac{5103}{2}\sqrt{3}s'+\frac{1701}{2}\sqrt{3}s-\frac{22411}{2}\sqrt{3}+\frac{53465}{4}\right)}{\pi z^8}\\
+\frac{1}{19208}\sqrt{3} \frac{\left(-\frac{413343}{64}s'+\frac{7490259}{64}s+\frac{76545}{64}\sqrt{3}s'+\frac{66339}{64}\sqrt{3}s+\frac{7226181}{128}\sqrt{3}-\frac{4527333}{32}\right)}{\pi z^8}\\
 - \frac{3}{19208} \frac{\left(-\frac{7582113}{64}s' + \frac{15309}{64}s+/frac{5103}{64}\sqrt{3}s'+\frac{15309}{64}\sqrt{3}s+\frac{295803}{128} \sqrt{3}+\frac{2697651}{32}\right)}{\pi z^8}
\end{smallmatrix}\right)\\
& +\quad O(\frac{1}{z^9})\,.
\end{array}
$$

\item The expansion of vector fields $\vecZ^1_5$ and $\vecZ_5^2$ are expressed by,

$$
\vecZ_5^1(\vecX_0) = \left(\begin{smallmatrix}
-\frac{64881}{1404928}\frac{\sqrt{3}}{z^4}+\frac{5000643}{157351936}\frac{\sqrt{3}}{z^5}-\frac{1007277825}{1258815488}\frac{\sqrt{3}}{z^6} -\frac{7453380147}{1258815488}\frac{\sqrt{3}}{z^7}-\frac{9398907099}{1258815488}\frac{\sqrt{3}}{z^8}\\
0\\
\frac{65}{2592}\sqrt{3}+\frac{21627}{5619712}\frac{\sqrt{3}}{z^4}-\frac{67204577}{629407744}\frac{\sqrt{3}}{z^5}+\frac{2318031745}{15105785856}\frac{\sqrt{3}}{z^6}-\frac{1808333293}{15105785856}\frac{\sqrt{3}}{z^7}-\frac{31796394461}{15105785856}\frac{\sqrt{3}}{z^8}\\
\frac{65}{864}+\frac{64881}{5619712 z^4} -\frac{201613731}{629407744 z^5}+\frac{2318031745}{5035261952 z^6}-/frac{1808333293}{5035261952 z^7}-\frac{31796394461}{5035261952 z^8}\\
\frac{65}{1296}-\frac{64881}{2809856 z^4}-\frac{40022909}{314703872 z^5}-\frac{2121362435}{7552892928 z^6}+\frac{1191556679}{7552892928 z^7}+\frac{62259224983}{7552892928 z^8}
\end{smallmatrix}\right)\,,
$$

%$$
%\vecZ_5^2(\vecX_0,s,s') = \left(\begin{smallmatrix}
%-\frac{6561}{10976}\frac{\left(\frac{1}{2}+\frac{1}{2}\sqrt{3}\right)\left(1-\frac{1}{4}\sqrt{3}\right)}{\pi z^7}+\frac{27}{38416}\frac{\left(\frac{1}{2}+\frac{1}{2}\sqrt{3}\right)\left(1-\frac{1}{4}\sqrt{3}\right)\left(-\frac{5103}{4}s'-\frac{1701}{4}s+\frac{2463}{64}\sqrt{3}+\frac{1431799}{64}\right)}{\pi z^8}\\
%\frac{6561}{10976} \frac{\left(-\frac{1}{2}\sqrt{3}+\frac{1}{2}\right)}{\pi z^7}-\frac{27}{38416} \frac{\left(-\frac{1701}{8}s'-\frac{1701}{8}s+\frac{5103}{8}\sqrt{3}s'-\frac{1701}{8}\sqrt{3}s-\frac{171863}{16}\sqrt{3}+\frac{657773}{64}\right)}{\pi z^8}\\
%-\frac{27}{21952}\frac{\left(\frac{729}{16}\sqrt{3}-\frac{1053}{16}\right)}{\pi z^7}+\frac{1}{19208}\sqrt{3} \frac{\left(-\frac{413343}{256}s'-\frac{28812861}{256}s+\frac{76545}{256}\sqrt{3}s'+\frac{66339}{256}\sqrt{3}s+\frac{102986487}{2048}\sqrt{3}+\frac{144659439}{1024}\right)}{\pi z^8}\\
% \frac{27}{21952}\sqrt{3} \frac{\left(\frac{567}{16}\sqrt{3}-\frac{243}{16}\right)}{\pi z^7} - \frac{3}{19208} \frac{\left(\frac{28721007}{256}s'+\frac{15309}{256}s+\frac{5103}{256}\sqrt{3}s'+\frac{15309}{256}\sqrt{3}s-\frac{7786791}{2048}\sqrt{3}-\frac{403566543}{1024}\right)}{\pi z^8}\\
%-\frac{2187}{21952}\frac{\left(\frac{1}{2}+\frac{1}{2}\sqrt{3}\right)\left(1-\frac{1}{4}\sqrt{3}\right)\sqrt{3}}{\pi z^7}+\frac{9}{76832}\frac{\left(\frac{1}{2}+\frac{1}{2}\sqrt{3}\right)\left(1-\frac{1}{4}\sqrt{3}\right)\sqrt{3} \left(-\frac{5103}{4}s'-\frac{1701}{4}s+\frac{2463}{64}\sqrt{3}+\frac{1431799}{64}\right)}{\pi z^8}
%\end{smallmatrix}\right)\,.
%$$

$$
\begin{array}{lll}
\vecZ_5^2(\vecX_0,s,s') =& \left(\begin{smallmatrix}
-\frac{6561}{10976}\frac{\left(\frac{1}{2}+\frac{1}{2}\sqrt{3}\right)\left(1-\frac{1}{4}\sqrt{3}\right)}{\pi z^7}\\
\frac{6561}{10976} \frac{\left(-\frac{1}{2}\sqrt{3}+\frac{1}{2}\right)}{\pi z^7}\\
-\frac{27}{21952}\frac{\left(\frac{729}{16}\sqrt{3}-\frac{1053}{16}\right)}{\pi z^7}\\
 \frac{27}{21952}\sqrt{3} \frac{\left(\frac{567}{16}\sqrt{3}-\frac{243}{16}\right)}{\pi z^7}\\
-\frac{2187}{21952}\frac{\left(\frac{1}{2}+\frac{1}{2}\sqrt{3}\right)\left(1-\frac{1}{4}\sqrt{3}\right)\sqrt{3}}{\pi z^7}
\end{smallmatrix}\right)\quad +\\
&\left(\begin{smallmatrix}
+\frac{27}{38416}\frac{\left(\frac{1}{2}+\frac{1}{2}\sqrt{3}\right)\left(1-\frac{1}{4}\sqrt{3}\right)\left(-\frac{5103}{4}s'-\frac{1701}{4}s+\frac{2463}{64}\sqrt{3}+\frac{1431799}{64}\right)}{\pi z^8}\\
-\frac{27}{38416} \frac{\left(-\frac{1701}{8}s'-\frac{1701}{8}s+\frac{5103}{8}\sqrt{3}s'-\frac{1701}{8}\sqrt{3}s-\frac{171863}{16}\sqrt{3}+\frac{657773}{64}\right)}{\pi z^8}\\
+\frac{1}{19208}\sqrt{3} \frac{\left(-\frac{413343}{256}s'-\frac{28812861}{256}s+\frac{76545}{256}\sqrt{3}s'+\frac{66339}{256}\sqrt{3}s+\frac{102986487}{2048}\sqrt{3}+\frac{144659439}{1024}\right)}{\pi z^8}\\
  - \frac{3}{19208} \frac{\left(\frac{28721007}{256}s'+\frac{15309}{256}s+\frac{5103}{256}\sqrt{3}s'+\frac{15309}{256}\sqrt{3}s-\frac{7786791}{2048}\sqrt{3}-\frac{403566543}{1024}\right)}{\pi z^8}\\
+\frac{9}{76832}\frac{\left(\frac{1}{2}+\frac{1}{2}\sqrt{3}\right)\left(1-\frac{1}{4}\sqrt{3}\right)\sqrt{3} \left(-\frac{5103}{4}s'-\frac{1701}{4}s+\frac{2463}{64}\sqrt{3}+\frac{1431799}{64}\right)}{\pi z^8}
\end{smallmatrix}\right)\\
& +\quad O(\frac{1}{z^9})\,.
\end{array}
$$

\end{itemize}
Then, we express the function $\det_{int}$ computed at $(\vecX_0,s,s',t,t')$

$$
\begin{array}{cc}
\det_{int}(\vecX_0,s,s',t,t') =& -\frac{4209544161}{5279854836580352} \left(-623289 s'+623289 t+3220141 t'\right.\\
&-3220141 s+1153029 s^2+384343s'^2\sqrt{3}\\
&+623289\sqrt{3}s'-384343 s^2-384343\sqrt{3}s'^2\\
&-623289\sqrt{3} t-1682769 \sqrt{3}t'+1682769\sqrt{3}s\\
&\left.+384343 t'^2\sqrt{3}+384343\sqrt{3} t^2-1153029 t'^2-384343 t^2\right)\frac{1}{\pi^2 z^{24}} \\
&+ O(\frac{1}{z^{25}})\,.
\end{array}
$$

This formal expressions lead to conclude the proof of theorem \ref{thm:three_sphere}.

\end{document}